\documentclass[11pt]{article}
\PassOptionsToPackage{obeyspaces}{url}
\usepackage[hidelinks]{hyperref}
\usepackage{bookmark}
\bookmarksetup{open,numbered,depth=5}
\usepackage{amsfonts}
\usepackage{amsmath}
\usepackage{amssymb, amsthm, enumitem}
\usepackage{breakurl}
\usepackage{pgf,tikz}
\usepackage{bm}
\usepackage{cancel}
\usepackage{comment}
\usetikzlibrary{calc}

\usepackage{etoolbox} 
\patchcmd{\thebibliography}{\leftmargin\labelwidth}{\leftmargin\labelwidth\addtolength\itemsep{-0.1\baselineskip}}{}{}

\usepackage{mathtools}
\usepackage{xstring}

\author{Zichao Dong\thanks{Department of Mathematical Sciences, Carnegie Mellon University, Pittsburgh, PA 15213, USA\@. Supported in part by U.S.\ taxpayers through NSF grant DMS-2154063. \texttt{zichaod@andrew.cmu.edu}. } \and Zijian Xu\thanks{School of Mathematical Sciences, Peking University, Beijing 100871. \texttt{2200010770@stu.pku.edu.cn}.}}

\title{Rainbow even cycles}
\date{}

\usepackage[nameinlink]{cleveref}

\newtheorem{theorem}{Theorem}
\newtheorem{lemma}[theorem]{Lemma}
\newtheorem{corollary}[theorem]{Corollary}
\newtheorem{observation}[theorem]{Observation}
\newtheorem{proposition}[theorem]{Proposition}

\newcommand*{\eqdef}{\stackrel{\mbox{\normalfont\tiny def}}{=}} 
\DeclarePairedDelimiter\abs{\lvert}{\rvert}                     
\DeclarePairedDelimiter\bracket{\langle}{\rangle}                     

\newcommand*{\N}{\mathbb{N}}                                    
\newcommand*{\Z}{\mathbb{Z}}                                    
\newcommand*{\R}{\mathbb{R}}                                    
\newcommand*{\cD}{\mathcal{D}}

\newcommand*{\sB}{\mathsf{B}}
\newcommand*{\sC}{\mathsf{C}}
\newcommand*{\sE}{\mathsf{E}}
\newcommand*{\sD}{\mathsf{D}}
\newcommand*{\fF}{\mathfrak{F}}
\newcommand*{\sF}{\mathsf{F}}
\newcommand*{\sP}{\mathsf{P}}
\newcommand*{\sQ}{\mathsf{Q}}
\newcommand*{\sT}{\mathsf{T}}

\newcommand*{\cP}{\mathcal{P}}
\newcommand*{\sS}{\mathsf{S}}
\newcommand*{\se}{\mathsf{e}}
\renewcommand*{\sf}{\mathsf{f}}
\newcommand*{\sg}{\mathsf{g}}
\newcommand*{\sG}{\mathsf{G}}
\newcommand*{\sh}{\mathsf{h}}

\renewcommand*{\sp}{\mathsf{p}}
\newcommand*{\sq}{\mathsf{q}}

\newcommand*{\sk}{\mathsf{k}}
\newcommand*{\cF}{\mathcal{F}}

\newcommand*{\sX}{\mathsf{X}}

\newcommand*{\fw}{Frankenstein }

\DeclareMathOperator{\depth}{depth}
\DeclareMathOperator{\Depth}{Depth}

\DeclareMathOperator{\child}{Child}

\DeclareMathOperator{\dist}{dist}

\newcommand{\lc}{\left\lceil}
\newcommand{\rc}{\right\rceil}
\newcommand{\lf}{\left\lfloor}
\newcommand{\rf}{\right\rfloor}

\crefname{enumi}{step}{steps}
\crefname{part}{part}{parts}

\allowdisplaybreaks[4]

\begin{document}
	
	\maketitle
	
	\begin{abstract}
		We prove that every family of (not necessarily distinct) even cycles $D_1, \dotsc, D_{\lf 1.2(n-1) \rf+1}$ on some fixed $n$-vertex set has a rainbow even cycle (that is, a set of edges from distinct $D_i$'s, forming an even cycle). This resolves an open problem of Aharoni, Briggs, Holzman and Jiang. Moreover, the result is best possible for every positive integer $n$. 
	\end{abstract}
	
	\section{Introduction} \label{sec:intro}
	
	Let $\cF$ be a set family. A rainbow set with respect to $\cF$ is a subset $R$ (without repeated elements) of $\cup \cF$ (i.e.~$\bigcup_{F \in \cF} F$) such that there exists an injection $\sigma \colon R \to \cF$ with $r \in \sigma(r)$ for all $r \in R$. In other words, each element $r \in R$ comes from a distinct $F \in \cF$. We think about each set in $\cF$ as a different color class, and hence use the term ``rainbow''. An important remark here is that a ``family'' refers to a ``multiset'', since an element in $\cup \mathcal{F}$ may appear with more than one colors. 
	
	Suppose every $F \in \cF$ satisfies property $\cP$. What is the minimum size of $\cF$ such that a rainbow subset of $\cup \cF$ satisfying $\cP$ always exists? One famous result of this type is the colorful version of Carath\'{e}odory's theorem due to B\'{a}r\'{a}ny \cite{barany}, which asserts that every family of $n+1$ subsets of $\R^n$, each containing a point $p$ in its convex hull, has a rainbow subset whose convex hull contains $p$ as well. Such problems are also studied in graph theory. Aharoni and Berger \cite{aharoni_berger} proved that any family of $2n - 1$ matchings of size $n$ in a bipartite graph contains a rainbow matching of size $n$. Other results of this type on cycles and triangles can be found in \cite{aharoni_briggs_holzman_jiang,gyori,goorevitch_holzman}. 
	
	There are studies of rainbow graphs in a different context: Given an edge-colored graph, what conditions guarantee a certain subgraph whose edges have distinct colors? Due to the relation with Latin squares, rainbow matchings have received extensive attention. See \cite{aharoni_berger_chudnovsky_zerbib,correia_pokrovskiy_sudakov} for recent works. As a starting point to find colorful variants of Tur\'{a}n's theorem, the existence of rainbow triangles is analyzed in \cite{aharoni_devos_gonzalez_montejano_samal,aharoni_devos_holzman}. A rainbow version of Dirac's theorem on Hamiltonian cycles can be found in \cite{joos_kim}. 
	
	\smallskip
	
	Throughout the paper, a graph, without further specification, refers to a simple graph $\sG$ which is a set of colored edges. Formally, $\sG$ is a set of pairs $\se = (e, \alpha)$ where $e$'s are distinct edges (i.e.~different pairs of two distinct vertices) and $\alpha$'s are (not necessarily distinct) colors. For $\se = (uv, \alpha) \in \sG$ where $uv \eqdef \{u, v\}$, denote $V(\se) \eqdef \{u, v\}, \, \chi(\se) = \alpha$. Then write $V(\sG) \eqdef \bigcup_{\se \in \sG} V(\se), \, E(\sG) \eqdef \{V(\se) : \se \in \sG\}$ and $\chi(\sG) \eqdef \{\chi(\se) : \se \in \sG\}$ for the vertex set, the (uncolored) edge set and the color set, respectively. 
	
	Two edges $\se_1, \se_2$ are \emph{coincident} if they are of different colors and are on the same vertex set. That is, $V(\se_1) = V(\se_2)$ yet $\chi(\se_1) \ne \chi(\se_2)$. For two graphs $\sG_1, \sG_2$, we call them \emph{coincident} if there exists a bijection $\varphi \colon \sG_1 \to \sG_2$ such that $\se$ is coincident to $\varphi(\se)$ for all $\se \in \sG_1$. Note that coincident edges do not exist in a graph, since graphs are assumed to be simple. 
	
	This paper is devoted to the existence of a rainbow even cycle in a family of even cycles. A cycle is a graph $\sC$ such that its edges $E(\sC)$, viewed as an uncolored simple graph, form a cycle. In other words, $\sC = \{(v_1v_2, \alpha_1), \dotsc, (v_{\ell-1}v_{\ell}, \alpha_{\ell-1}), (v_{\ell}v_1, \alpha_{\ell})\}$ where $v_1, \dotsc, v_{\ell}$ are distinct and $\ell \ge 3$ is called the \emph{length} of $\sC$. For any $A \subseteq \{3, 4, 5, \dotsc\}$, an $A$-\emph{cycle} is a cycle whose length is some number from $A$. For example, an \emph{odd cycle}, a cycle of odd length, is a $\{3, 5, 7, \dotsc\}$-cycle. Similarly, an \emph{even cycle}, a cycle of even length, is a $\{4, 6, 8, \dotsc\}$-cycle. For any integer $k \ge 3$,  a $k$-cycle refers to a $\{k\}$-cycle. 
	
	Hereafter a family $\cF = \{\sE_1, \dotsc, \sE_m\}$ is a family of cycles. We remark that $\cF$ is a family implicitly implies that $\chi(\sE_i) = \{\alpha_i\}$ while $\alpha_1, \dotsc, \alpha_m$ are distinct. Since each $\sE_i$ is a monochromatic cycle, we view $\cF$ as an edge-colored multi-graph (i.e.~a set of colored edges where coincident edges are allowed). A subgraph of $\cF$ is then a graph $\sE$ where $\sE \subseteq \bigcup_{i=1}^m \sE_i$. In \hyperlink{figone}{Figure~1}, the family $\cD = \{\sD_1, \sD_2, \sD_3, \sD_4\}$ consists of four $4$-cycles on seven vertices where $\sD_2, \sD_3$ are coincident. Let $\chi(\sD_i) = \alpha_i \, (i = 1, 2, 3, 4)$. Then $\sD \eqdef \{(v_0v_1, \alpha_1), (v_1v_2, \alpha_2), (v_2v_3, \alpha_3), (v_3v_0, \alpha_4)\}$ is a rainbow $4$-cycle subgraph of $\cD$. 
	
	\begin{center}
		\begin{tikzpicture}[x=2cm,y=1.2cm, scale = 1]
			\clip(-4,-1.9) rectangle (4,1.4);
			\draw [color=red] (-0.6,-1) -- (-1.2,0) -- (-0.6,1) -- (0,0) -- (-0.6,-1);
			\draw [color=blue] (0.6,-1) -- (1.2,0) -- (0.6,1) -- (0,0) -- (0.6,-1);
			\draw [color=green] (-0.625,-1.05) -- (-0.625,1.05) -- (0.625,1.05) -- (0.625,-1.05) -- (-0.625,-1.05);
			\draw [color=violet] (-0.575,-0.95) -- (-0.575,0.95) -- (0.575,0.95) -- (0.575,-0.95) -- (-0.575,-0.95);
			\draw [fill=black] (-1.2,0) circle (4pt);
			\draw [fill=black] (-0.6,-1) circle (4pt);
			\draw [fill=black] (-0.6,1) circle (4pt);
			\draw [fill=black] (0,0) circle (4pt);
			\draw [fill=black] (0.6,-1) circle (4pt);
			\draw [fill=black] (0.6,1) circle (4pt);
			\draw [fill=black] (1.2,0) circle (4pt);
			\node at (0,-0.3) {$v_0$};
			\node at (-0.75,1.15) {$v_1$};
			\node at (-0.75,-1.15) {$v_2$};
			\node at (0.75,-1.15) {$v_3$};
			\node at (0.75,1.15) {$v_4$};
			\node at (-1.375,0) {$v_5$};
			\node at (1.375,0) {$v_6$};
			\begin{footnotesize}
				\node at (-1,0.65) {$\textcolor{red}{\sD_1}$};
				\node at (0,0.775) {$\textcolor{violet}{\sD_2}$};
				\node at (0,1.225) {$\textcolor{green}{\sD_3}$};
				\node at (1,0.65) {$\textcolor{blue}{\sD_4}$};
			\end{footnotesize}
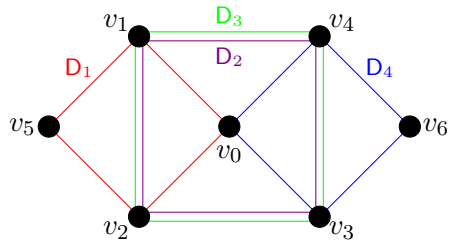
			\node at (0,-1.7) {\textbf{\hypertarget{figone}{Figure 1:}} An example family $\cD$ viewed as an edge-colored multigraph. };
		\end{tikzpicture}
	\end{center}
	
	We shall say that a family $\mathcal{F}$ \emph{contains} a graph $\sG$ if $\sG$ is a subgraph of $\cF$. 
	
	\begin{theorem} \label{thm:rainbow_oc}
		\emph{(\cite{aharoni_briggs_holzman_jiang})} Every family of $2\lc\frac{n}{2}\rc-1$ odd cycles on $n$ vertices contains a rainbow odd cycle. 
	\end{theorem}
	
	The tightness of \Cref{thm:rainbow_oc} is witnessed by a family of $2\bigl( \bigl\lceil \frac{n}{2} \bigr\rceil - 1 \bigr)$ many coincident odd cycles on $2 \bigl\lceil \frac{n}{2} \bigr\rceil - 1$ vertices. As for even cycles, Aharoni, Briggs, Holzman and Jiang also deduced in \cite{aharoni_briggs_holzman_jiang} that the maximum size of a family on $n$ vertices containing no rainbow even cycle is between roughly $\frac{6}{5}n$ and $\frac{3}{2}n$, and left the determination of the exact extremal number as an open problem. We answer this question by proving the following result: 
	
	\begin{theorem} \label{thm:rainbow_ec}
		Every family of $\bigl\lfloor \frac{6(n-1)}{5} \bigr\rfloor + 1$ even cycles on $n$ vertices contains a rainbow even cycle. 
	\end{theorem}
	
	The tightness of \Cref{thm:rainbow_ec} for each $n \ge 4$ (no even cycle exists when $n \le 3$) is seen as follows: The families $\cD_4, \cD_5, \cD_6, \cD_7, \cD_8$ in \hyperlink{figtwo}{Figure~2} are tight examples for $n = 4, 5, 6, 7, 8$, respectively. For larger $n$, we observe that by gluing together $\cD_{n-5}$ (a tight example for $n-5$) and $\cD_6$ at exactly one vertex (edge-disjoint henceforth) the resulting family $\cD_n$ is tight for $n$. We remark that the family $\cD_6$ and the inductive argument are already presented in \cite{aharoni_briggs_holzman_jiang}. 
	
	\begin{center}
		\begin{tikzpicture}[x=1.0cm,y=1cm]
			\clip(-1.5,-2) rectangle (1.5,1.5);
			\draw [color=red] (-1,-1)--(-1,1);
			\draw [color=red] (-1.1,-1)--(-1.1,1);
			\draw [color=red] (-0.9,-1)--(-0.9,1);
			\draw [color=red] (-1,1)--(1,1);
			\draw [color=red] (-1,0.9)--(1,0.9);
			\draw [color=red] (-1,1.1)--(1,1.1);
			\draw [color=red] (1,1)--(1,-1);
			\draw [color=red] (1.1,1)--(1.1,-1);
			\draw [color=red] (0.9,1)--(0.9,-1);
			\draw [color=red] (1,-1)--(-1,-1);
			\draw [color=red] (1,-1.1)--(-1,-1.1);
			\draw [color=red] (1,-0.9)--(-1,-0.9);
			\draw [fill=black] (-1,-1) circle (5pt);
			\draw [fill=black] (-1,1) circle (5pt);
			\draw [fill=black] (1,-1) circle (5pt);
			\draw [fill=black] (1,1) circle (5pt);
			\draw [color=black] (0,-1.7) node {$\cD_4$};
		\end{tikzpicture}
		\begin{tikzpicture}[x=1.0cm,y=1cm]
			\clip(-1.5,-2) rectangle (1.5,1.5);
			\draw [color=red] (-1.2,-0.95)--(0,-0.95);
			\draw [color=red] (-1.2,-1.05)--(0,-1.05);
			\draw [color=blue] (0,-0.95)--(1.2,-0.95);
			\draw [color=blue] (0,-1.05)--(1.2,-1.05);
			\draw [color=red] (-0.6,1.12)--(0.6,1.12);
			\draw [color=red] (-0.6,1.04)--(0.6,1.04);
			\draw [color=blue] (-0.6,0.96)--(0.6,0.96);
			\draw [color=blue] (-0.6,0.88)--(0.6,0.88);
			\draw [color=red] (-1.255,-1)--(-0.655,1);
			\draw [color=red] (-1.145,-1)--(-0.545,1);
			\draw [color=red] (-0.055,-1)--(0.545,1);
			\draw [color=red] (0.055,-1)--(0.655,1);
			\draw [color=blue] (1.255,-1)--(0.655,1);
			\draw [color=blue] (1.145,-1)--(0.545,1);
			\draw [color=blue] (-0.055,-1)--(-0.655,1);
			\draw [color=blue] (0.055,-1)--(-0.545,1);
			\draw [fill=black] (-1.2,-1) circle (5pt);
			\draw [fill=black] (0,-1) circle (5pt);
			\draw [fill=black] (1.2,-1) circle (5pt);
			\draw [fill=black] (-0.6,1) circle (5pt);
			\draw [fill=black] (0.6,1) circle (5pt);
			\draw [color=black] (0,-1.7) node {$\cD_5$};
		\end{tikzpicture}
		\begin{tikzpicture}[x=1.0cm,y=1cm]
			\clip(-1.5,-2) rectangle (1.5,1.5);
			\foreach \x in {-0.1,0,0.1}
			\foreach \y in {-1,1}
			\draw [color=red] (-1.2,\x+\y)--(0,\x+\y);
			\foreach \x in {-0.12,0,0.12}
			\draw [color=red] (-1.2+\x,-1)--(\x,1);
			\foreach \x in {-0.12,0,0.12}
			\draw [color=red] (-1.2+\x,1)--(\x,-1);
			\foreach \x in {-0.1,0,0.1}
			\foreach \y in {0,1.2}
			\draw [color=blue] (\x+\y,-1)--(\x+\y,1);
			\foreach \x in {-0.12,0,0.12}
			\draw [color=blue] (1.2+\x,-1)--(\x,1);
			\foreach \x in {-0.12,0,0.12}
			\draw [color=blue] (1.2+\x,1)--(\x,-1);
			\draw [fill=black] (-1.2,-1) circle (5pt);
			\draw [fill=black] (0,-1) circle (5pt);
			\draw [fill=black] (1.2,-1) circle (5pt);
			\draw [fill=black] (-1.2,1) circle (5pt);
			\draw [fill=black] (0,1) circle (5pt);
			\draw [fill=black] (1.2,1) circle (5pt);
			\draw [color=black] (0,-1.7) node {$\cD_6$};
		\end{tikzpicture}
		\begin{tikzpicture}[x=1.0cm,y=1cm]
			\clip(-1.5,-2) rectangle (1.5,1.5);
			\foreach \x in {-0.12,0,0.12}
			\foreach \y in {-1.2,0}
			\foreach \z in {-1,1}
			\draw [color=red] (\x+\y,0)--(-0.6+\x,\z);
			\foreach \x in {-0.12,0,0.12}
			\foreach \y in {1.2,0}
			\foreach \z in {-1,1}
			\draw [color=blue] (\x+\y,0)--(0.6+\x,\z);
			\draw [color=green] (-0.6,-1)--(-0.6,1);
			\draw [color=green] (-0.6,1)--(0.6,1);
			\draw [color=green] (0.6,1)--(0.6,-1);
			\draw [color=green] (0.6,-1)--(-0.6,-1);
			\draw [fill=black] (-1.2,0) circle (5pt);
			\draw [fill=black] (-0.6,-1) circle (5pt);
			\draw [fill=black] (-0.6,1) circle (5pt);
			\draw [fill=black] (0,0) circle (5pt);
			\draw [fill=black] (0.6,-1) circle (5pt);
			\draw [fill=black] (0.6,1) circle (5pt);
			\draw [fill=black] (1.2,0) circle (5pt);
			\draw [color=black] (0,-1.7) node {$\cD_7$};
		\end{tikzpicture}
		\begin{tikzpicture}[x=1.0cm,y=1cm]
			\clip(-1.5,-2) rectangle (1.5,1.5);
			\foreach \x in {-0.1,-0.05,0,0.05,0.1}
			\foreach \y in {-1.2,-0.4}
			\draw [color=red] (\x+\y,-0.3)--(\x+\y,0.3);
			\foreach \x in {-0.11,-0.055,0,0.055,0.11}
			\draw [color=red] (-1.2+\x,0.3)--(-0.8+\x,1);
			\foreach \x in {-0.11,-0.055,0,0.055,0.11}
			\draw [color=red] (-0.4+\x,0.3)--(-0.8+\x,1);
			\foreach \x in {-0.11,-0.055,0,0.055,0.11}
			\draw [color=red] (-1.2+\x,-0.3)--(-0.8+\x,-1);
			\foreach \x in {-0.11,-0.055,0,0.055,0.11}
			\draw [color=red] (-0.4+\x,-0.3)--(-0.8+\x,-1);
			\foreach \x in {-0.07,0,0.07}
			\foreach \y in {-1,1}
			\draw [color=blue] (-0.8+\x,\y)--(1.2+\x,-\y);
			\foreach \x in {-0.05,0,0.05}
			\foreach \y in {-1,1}
			\draw [color=blue] (-0.8,\x+\y)--(1.2,\x+\y);
			\draw [fill=black] (-1.2,-0.3) circle (5pt);
			\draw [fill=black] (-1.2,0.3) circle (5pt);
			\draw [fill=black] (-0.8,-1) circle (5pt);
			\draw [fill=black] (-0.8,1) circle (5pt);
			\draw [fill=black] (-0.4,-0.3) circle (5pt);
			\draw [fill=black] (-0.4,0.3) circle (5pt);
			\draw [fill=black] (1.2,-1) circle (5pt);
			\draw [fill=black] (1.2,1) circle (5pt);
			\draw [color=black] (0,-1.7) node {$\cD_8$};
		\end{tikzpicture}
		\begin{tikzpicture}
			\node at (0, 0) {\textbf{\hypertarget{figtwo}{Figure 2:}} Tight examples of \Cref{thm:rainbow_ec} for small $n$. };
		\end{tikzpicture}
	\end{center}
	
	\paragraph{Proof strategy.} To explain the strategy of our proof, we begin with a baby version of \Cref{thm:rainbow_ec} whose tightness is witnessed by, for example, a family of $n-1$ coincident Hamiltonian cycles. 
	
	\begin{theorem} \label{thm:rainbow_c}
		\emph{(\cite[Proposition~3.2]{aharoni_briggs_holzman_jiang})} Every family of $n$ cycles on $n$ vertices contains a rainbow cycle. 
	\end{theorem}
	
	\begin{proof}
		Let $\cF$ be such a family and $\sF$ be a maximal rainbow forest subgraph of $\cF$. Then $|\sF| \le n-1$, and so there is another edge $\se$, not coincident to any edge of $\sF$, whose color does not appear in $\sF$. The maximality of $\sF$ implies that $\se$ completes a rainbow cycle in the graph $\sF \cup \{\se\}$. 
	\end{proof}
	
	All these proofs proceed by first finding a \emph{spanning structure} $S$ (the rainbow forest $\sF$ in the proof above) and then analyzing another edge with an absent color in $S$. The proof of \Cref{thm:rainbow_oc} also uses a maximal rainbow forest as $S$. However, to prove \Cref{thm:rainbow_ec} we need some new spanning structure. 
	
	It turns out that $5$-cycles play a central role in the $\frac{6}{5}n$ upper bound. We thus call a cycle \emph{long} if its length is at least $6$. In particular, a rainbow $\{7, 9, 11, \dotsc\}$-cycle is a long rainbow odd cycle. Then our spanning structure, which we call \emph{Frankenstein graphs}, are (informally speaking) obtained by recursively at single vertices gluing together a collection of long rainbow odd cycles, rainbow trees, and another class of graphs named \emph{bad pieces}. 
	
	We shall formally define and characterize bad pieces and Frankenstein graphs in \Cref{sec:fgraph}. Then \Cref{sec:proof_ec} is devoted to the proof of \Cref{thm:rainbow_ec}. 
	
	\section{\fw graphs} \label{sec:fgraph}
	
	A path graph of length $k$ is a graph of the form
	\[
	\sP = \{(v_0v_1, \alpha_1), (v_1v_2, \alpha_2), \dotsc, (v_{k-1}v_k, \alpha_k)\}
	\]
	where $v_0, \dotsc, v_k$ are distinct. A \emph{theta graph} is a union of $3$ paths that share exactly their terminals. Formally, $\sG$ is a theta graph if $\sG = \sP_1 \cup \sP_2 \cup \sP_3$ where $\sP_1, \sP_2, \sP_3$ are paths with terminals $s, t$ and
	\begin{align*}
		V(\sP_1) \cap V(\sP_2) = V(\sP_2) \cap V(\sP_3) &= V(\sP_3) \cap V(\sP_1) = \{s, t\}, \\
		E(\sP_1) \cap E(\sP_2) = E(\sP_2) \cap E(\sP_3) &= E(\sP_3) \cap E(\sP_1) = \varnothing. 
	\end{align*}
	We use the name ``theta'' because one natural drawing of such a graph looks exactly like the Greek letter $\Theta$. See \hyperlink{figthree}{Figure~3} below as an illustration. 
	
	\vspace{-0.5em}
	\begin{center}
		\begin{tikzpicture}
			[x=0.8cm,y=0.8cm]
			\clip(-5,-2.05) rectangle (11,1.5);
			\draw [fill=black] (0,0) circle (1pt);
			\draw [fill=black] (1.5,0) circle (1pt);
			\draw [fill=black] (3,0) circle (1pt);
			\draw [fill=black] (4.5,0) circle (1pt);
			\draw [fill=black] (6,0) circle (1pt);
			\draw [fill=black] (2,1) circle (1pt);
			\draw [fill=black] (4,1) circle (1pt);
			\draw [fill=black] (1.2,-1) circle (1pt);
			\draw [fill=black] (2.4,-1) circle (1pt);
			\draw [fill=black] (3.6,-1) circle (1pt);
			\draw [fill=black] (4.8,-1) circle (1pt);
			\draw (0,0) -- (2,1);
			\draw (2,1) -- (4,1);
			\draw (4,1) -- (6,0);
			\draw (6,0) -- (4.8,-1);
			\draw (4.8,-1) -- (3.6,-1);
			\draw (3.6,-1) -- (2.4,-1);
			\draw (2.4,-1) -- (1.2,-1);
			\draw (1.2,-1) -- (0,0);
			\draw (0,0) -- (1.5,0);
			\draw (1.5,0) -- (3,0);
			\draw (3,0) -- (4.5,0);
			\draw (4.5,0) -- (6,0);
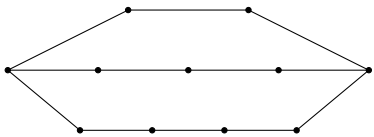
			\node at (3, -1.8) {\textbf{\hypertarget{figthree}{Figure 3:}} A theta graph on paths of lengths $3,4,5$, respectively. };
		\end{tikzpicture}
	\end{center}
	
	\begin{observation} \label{obs:theta_ec}
		Every rainbow theta graph has a rainbow even cycle subgraph. 
	\end{observation}
	
	\begin{proof}
		Suppose $\sP_1 \cup \sP_2 \cup \sP_3$ is a theta graph where $\sP_1, \sP_2, \sP_3$ are paths of common terminals. Then two of the paths, say $\sP_1$ and $\sP_2$, have lengths of the same parity, and so $\sP_1 \cup \sP_2$ is an even cycle. 
	\end{proof}
	
	We call a graph $\sG$ \emph{almost rainbow} if $|\chi(\sG)| = |\sG|-1$. That is, exactly two edges receive a same color, and the color of every other edge is unique. We call $\sB$ a \emph{bad piece} if $\sB$ is an almost rainbow theta graph on $3$ rainbow paths (sharing terminals) such that $|V(\sB)| \ge 6$. 
	
	\begin{center}
		\begin{tikzpicture}[x=1.0cm,y=0.8cm]
			\clip(-2.5,-2) rectangle (2.5,4.5);
			\draw  (-1.8,1.5)-- (0,4);
			\draw  (0,4)-- (1.8,1.5);
			\draw  (1.8,1.5)-- (0,-1);
			\draw  (0,-1)-- (-1.8,1.5);
			\draw  (-1.8,1.5)-- (1.8,1.5);
			\node at (0, -1.7) {\textbf{\hypertarget{figfourA}{Figure 4.A}}};
			\draw [fill=black] (-1.8,1.5) circle (1.0pt);
			\draw[color=black] (-2.05,1.5) node {$v_1$};
			\draw [fill=black] (1.8,1.5) circle (1.0pt);
			\draw[color=black] (2.05,1.5) node {$v_3$};
			\draw [fill=black] (0,-1) circle (1.0pt);
			\draw[color=black] (0,-1.21) node {$v_2$};
			\draw [fill=black] (0,4) circle (1.0pt);
			\draw[color=black] (0,4.2) node {$v_4$};
			\draw[color=black] (0,1.7) node {$a$};
			\draw[color=black] (0.8,0.4) node {$c$};
			\draw[color=black] (-0.8,0.45) node {$b$};
			\draw[color=black] (-0.825,2.6) node {$c$};
			\draw[color=black] (0.75,2.6) node {$d$};
		\end{tikzpicture}
		\qquad
		\begin{tikzpicture}[x=1.0cm,y=0.8cm]
			\clip(-2.5,-2) rectangle (2.5,4.5);
			\draw  (0.,2.)-- (-2.,1.);
			\draw  (-2.,1.)-- (-1.,-1.);
			\draw  (-1.,-1.)-- (1.,-1.);
			\draw  (1.,-1.)-- (2.,1.);
			\draw  (2.,1.)-- (0.,2.);
			\draw  (0.,4.)-- (0.,2.);
			\draw  (0.,4.)-- (2.,1.);
			\draw [fill=black] (-1.,-1.) circle (1.0pt);
			\draw[color=black] (-0.8,-0.8) node {$v_3$};
			\draw [fill=black] (1.,-1.) circle (1.0pt);
			\draw[color=black] (0.8,-0.8) node {$v_4$};
			\draw [fill=black] (-2.,1.) circle (1.0pt);
			\draw[color=black] (-1.7,0.9) node {$v_2$};
			\draw [fill=black] (2.,1.) circle (1.0pt);
			\draw[color=black] (1.7,0.9) node {$v_5$};
			\draw [fill=black] (0.,2.) circle (1.0pt);
			\draw[color=black] (0,1.75) node {$v_1$};
			\draw[color=black] (-0.9,1.35) node {$a$};
			\draw[color=black] (-1.325,0.075) node {$b$};
			\draw[color=black] (0,-0.8) node {$c$};
			\draw[color=black] (1.325,0.075) node {$d$};
			\draw[color=black] (0.9, 1.35) node {$e$};
			\draw [fill=black] (0.,4.) circle (1.0pt);
			\draw[color=black] (0.,4.2) node {$v_6$};
			\draw[color=black] (-0.175,3.) node {$a$};
			\draw[color=black] (1.175,2.65) node {$g$};
			\node at (0, -1.7) {\textbf{\hypertarget{figfourB}{Figure 4.B}}};
		\end{tikzpicture}
		\qquad
		\begin{tikzpicture}[x=1.0cm,y=0.8cm]
			\clip(-2.5,-2) rectangle (2.5,4.5);
			\draw  (0.,2.)-- (-2.,1.);
			\draw  (-2.,1.)-- (-1.,-1.);
			\draw  (-1.,-1.)-- (1.,-1.);
			\draw  (1.,-1.)-- (2.,1.);
			\draw  (2.,1.)-- (0.,2.);
			\draw  (0.,4.)-- (0.,2.);
			\draw  (0.,4.)-- (2.,1.);
			\draw [fill=black] (-1.,-1.) circle (1.0pt);
			\draw[color=black] (-0.8,-0.8) node {$v_3$};
			\draw [fill=black] (1.,-1.) circle (1.0pt);
			\draw[color=black] (0.8,-0.8) node {$v_4$};
			\draw [fill=black] (-2.,1.) circle (1.0pt);
			\draw[color=black] (-1.7,0.9) node {$v_2$};
			\draw [fill=black] (2.,1.) circle (1.0pt);
			\draw[color=black] (1.7,0.9) node {$v_5$};
			\draw [fill=black] (0.,2.) circle (1.0pt);
			\draw[color=black] (0,1.75) node {$v_1$};
			\draw[color=black] (-0.9,1.35) node {$a$};
			\draw[color=black] (-1.325,0.075) node {$b$};
			\draw[color=black] (0,-0.8) node {$c$};
			\draw[color=black] (1.35,0.075) node {$b$};
			\draw[color=black] (0.9, 1.35) node {$e$};
			\draw [fill=black] (0.,4.) circle (1.0pt);
			\draw[color=black] (0.,4.2) node {$v_6$};
			\draw[color=black] (-0.175,3.) node {$f$};
			\draw[color=black] (1.175,2.65) node {$g$};
			\node at (0, -1.7) {\textbf{\hypertarget{figfourC}{Figure 4.C}}};
		\end{tikzpicture}
	\end{center}
	
	For example, \hyperlink{figfourA}{Figure~4.A} is not a bad piece because it contains only $4$ vertices; \hyperlink{figfourB}{Figure~4.B} is a bad piece on $6$ vertices and $7$ edges consisting of rainbow paths $v_1v_5$, $v_1v_2v_3v_4v_5$ and $v_1v_6v_5$; \hyperlink{figfourC}{Figure~4.C} is not a bad piece because $v_1v_2v_3v_4v_5$ is not rainbow (as witnessed by $(v_2v_3, b)$ and $(v_4v_5, b)$). 
	
	\begin{observation} \label{obs:bestimate}
		If $\sB$ is a bad piece, then $|V(\sB)| = |\chi(\sB)| \le \frac{6}{5}(|V(\sB)|-1)$. 
	\end{observation}
	
	\begin{proof}
		Since $\sB$ is a theta graph, we have $|V(\sB)| = |\sB|-1$. Notice that $\sB$ is almost rainbow, we see that $|\chi(\sB)| = |\sB|-1$. It follows that $n \eqdef |\chi(\sB)| = |V(\sB)| \ge 6$, and hence $\frac{|\chi(\sB)|}{|V(\sB)|-1} = \frac{n}{n-1} \le \frac{6}{5}$. 
	\end{proof}
	
	\begin{observation} \label{obs:bflip}
		If $\sB$ is a bad piece, then for any distinct $v_1, v_2 \in V(\sB)$, there exists in $\sB$ a rainbow path subgraph whose terminals are $v_1$ and $v_2$. 
	\end{observation}
	
	\begin{proof}
		Since $|\chi(\sB)| = |\sB|-1$, it suffices to show that $v_1, v_2$ are vertices of a cycle in $\sB$. Suppose $\sB$ consists of three rainbow paths $\sP_1, \sP_2, \sP_3$. If $v_1$ and $v_2$ are on a same path, say $\sP_1$, then $\sP_1 \cup \sP_2$ is such a cycle. If $v_1$ and $v_2$ are on different paths, say $\sP_1$ and $\sP_2$, then $\sP_1 \cup \sP_2$ is such a cycle. 
	\end{proof}
	
	Let $G$ be a graph. We call $\cP = \{\sG_1, \dotsc, \sG_m\}$ a \emph{partition} if $\sG = \bigcup_{i=1}^m \sG_i$ and $|V(\sG_i) \cap V(\sG_j)| \le 1$, $\chi(\sG_i) \cap \chi(\sG_j) = \varnothing$ for every distinct $\sG_i, \sG_j$. We shall often abuse the notation by writing $\cP(\sG) = \cP$. Indeed, $\cP(\sG)$ is not a function of $\sG$, as the partition is usually not unique. The notation emphasizes that the partition is of $\sG$. In this sense, $\fF$ is a \emph{Frankenstein graph} if it admits a partition
	\[
	\cP(\fF) = \{\sC_1, \dotsc, \sC_c, \sB_1, \dotsc, \sB_b, \sT_1, \dotsc, \sT_t\} \qquad (c \ge 0, \, b \ge 0, \, t \ge 0, \, c+b+t \ge 1)
	\]
	where $\sC$'s are long rainbow odd cycles, $\sB$'s are bad pieces, and $\sT$'s are rainbow trees, such that
	\begin{enumerate}[label=(F\arabic*), ref=(F\arabic*)]
		\item\label{Fr_tree} $V(\sT_p) \cap V(\sT_q) = \varnothing$ for any distinct $p, q$, and
		\item\label{Fr_ecfree} no rainbow even cycle subgraph exists in $\fF$. 
	\end{enumerate}
	
	\begin{theorem} \label{thm:fsteps}
		For any Frankenstein graph $\fF$ with $\cP(\fF) = \{\sG_1, \dotsc, \sG_m\}$, there exists a permutation $\sigma$ on $[m]$ such that $\fF_i \eqdef \sG_{\sigma(1)} \cup \dotsb \cup \sG_{\sigma(i)}$ satisfies $|V(\fF_i) \cap V(\sG_{\sigma(i+1)})| \le 1$ for each $i \in [m-1]$. 
	\end{theorem}
	
	\Cref{thm:fsteps} suggests the following way to think about a connected Frankenstein graph $\fF$: Suppose the partition of $\fF$ is $\cP(\fF) = \{\sG_1, \dotsc, \sG_m\}$. Then one can order the parts as $\fF_1 \eqdef \sG_1', \sG_2', \dotsc, \sG_m'$, and recursively glue together $\sG_{i+1}'$ and the $i$'th graph $\fF_i$ at some single vertex to make the $(i+1)$'st graph $\fF_{i+1}$, such that eventually $\fF_m$ is exactly $\fF$. To prove \Cref{thm:fsteps}, we need some preparations. 
	
	\begin{lemma} \label{lem:ectest}
		Let $\sC$ be a rainbow cycle. Assume $\sX$ is a rainbow cycle or a bad piece with $\sX, \sC \setminus \sX$ being color-disjoint and $E(\sC) \setminus E(\sX) \ne \varnothing$. If $|V(\sC) \cap V(\sX)| \ge 2$, then $\sC \cup \sX$ contains a rainbow even cycle. 
	\end{lemma}
	
	Informally speaking, this technical result is helpful because it tells us that a rainbow cycle is likely to form a rainbow even cycle together with a long rainbow odd cycle or a bad piece. 
	
	\begin{proof} 
		Since $E(\sC) \setminus E(\sX) \ne \varnothing$, there exists an edge $\se \in \sC$ that is not coincident to any edge of $\sX$. Starting from $\se$ and moving along $\sC$ in opposite directions, we define the first vertices to meet on $\sX$ as $s_0, t_0$, thanks to $|V(\sC) \cap V(\sX)| \ge 2$. Then there exists a subpath $\sP_0$ (i.e.~a path subgraph) of $\sC \setminus \sX$ satisfying $\se \in \sP_0$. Here $s_0, t_0$ are terminals of $\sP_0$, $V(\sP_0) \cap V(\sX) = \{s_0, t_0\}$ and $\chi(\sP_0) \cap \chi(\sX) = \varnothing$. 
		
		We claim the existence of a rainbow theta subgraph in $\sX \cup \sP_0$, and so \Cref{obs:theta_ec} guarantees a rainbow even cycle subgraph in $\sX \cup \sC$. 
		
		If $\sX$ is a rainbow cycle, then $\sX \cup \sP_0$ is a rainbow theta graph. 
		
		If $\sX$ is a bad piece which consists of rainbow paths $\sP_1, \sP_2, \sP_3$ that share terminals $s$ and $t$, then $\sX \cup \sP_0$ is almost rainbow. In fact, we can always remove a subpath containing one of the repeated-\linebreak color edges on one of $\sP_1, \sP_2, \sP_3$ to get a rainbow theta graph. To be more specific, we assume without loss that the repeated color happens on $\sP_1$ and $\sP_3$. If $x,y \in V(\sP_i)$ for some fixed $i \in [3]$, then there exists a unique subpath of $\sP_i$ with terminals $x$ and $y$, and we denote by $\sP_{x, y}$ this subpath. 
		\begin{itemize}
			\item If $s_0$ and $t_0$ lie on a same $\sP_i$, then one of $V(\sP_1) \setminus \{s, t\}$ and $V(\sP_3) \setminus \{s, t\}$ is disjoint from $V(\sP_0)$, say $V(\sP_1) \setminus \{s, t\}$. This implies that $(\sP_2 \cup \sP_3) \cup \sP_0 \subseteq \sX \cup \sC$ is a rainbow theta graph. 
			\item Otherwise, at least one of $s_0$ and $t_0$ lies on $\sP_1 \cup \sP_3$, say $s_0 \in V(\sP_1)$. We further assume that the repeated-color edge (denoted by $*$) appears on $\sP_{s, s_0}$ rather than $\sP_{t, s_0}$ in $\sP_1$. See \hyperlink{figfive}{Figure~5}. 
			\begin{itemize}
				\item If $t_0 \in V(\sP_2)$, then by removing $\sP_3$ from $\sX \cup \sP_0$ we are left with a rainbow theta graph. 
				\item If $t_0 \in V(\sP_3)$, then by removing $\sP_{s, s_0}$ from $\sX \cup \sP_0$ we are left with a rainbow theta graph. 
			\end{itemize}
		\end{itemize}
		
		\begin{center}
			\begin{tikzpicture}
				[x=1.0cm,y=1.0cm]
				\clip(-0.5,-1.6) rectangle (6.5,1.6);
				\draw [color=blue] (1.5,0) -- (1.5,0.5) -- (2.2,0.5) -- (2.9,0.5) -- (3.6,0.5) -- (3.6,1);
				\draw [fill=black] (0,0) circle (1pt);
				\draw [fill=black] (1.5,0) circle (1pt);
				\draw [fill=black] (3,0) circle (1pt);
				\draw [fill=black] (4.5,0) circle (1pt);
				\draw [fill=black] (6,0) circle (1pt);
				\draw [fill=black] (1.2,1) circle (1pt);
				\draw [fill=black] (2.4,1) circle (1pt);
				\draw [fill=black] (3.6,1) circle (1pt);
				\draw [fill=black] (4.8,1) circle (1pt);
				\draw [fill=black] (1.2,-1) circle (1pt);
				\draw [fill=black] (2.4,-1) circle (1pt);
				\draw [fill=black] (3.6,-1) circle (1pt);
				\draw [fill=black] (4.8,-1) circle (1pt);
				\draw [fill=black] (1.5,0.5) circle (1pt);
				\draw [fill=black] (2.2,0.5) circle (1pt);
				\draw [fill=black] (2.9,0.5) circle (1pt);
				\draw [fill=black] (3.6,0.5) circle (1pt);
				\draw (0,0) -- (1.2,1) -- (2.4,1) -- (3.6,1) -- (4.8,1) -- (6,0);
				\draw[dash pattern=on 3pt off 3pt] (0,0) -- (1.2,-1) -- (2.4,-1) -- (3.6,-1) -- (4.8,-1) -- (6,0);
				\draw (0,0) -- (1.5,0);
				\draw (1.5,0) -- (3,0);
				\draw (3,0) -- (4.5,0);
				\draw (4.5,0) -- (6,0);
				\node at (-0.2,0) {$s$};
				\node at (6.2,0) {$t$};
				\node at (3.6,1.2) {$s_0$};
				\node at (1.5,-0.2) {$t_0$};
				\node at (2.7,1.225) {$\sP_1$};
				\node at (2.7,-0.2) {$\sP_2$};
				\node at (2.7,-1.2) {$\sP_3$};
				\begin{scriptsize}
					\node at (2.575,0.65) {\textcolor{blue}{$\sP_0$}};
				\end{scriptsize}
				\node at (1.8,1.15) {$*$};
				\node at (4.2,-0.85) {$*$};
			\end{tikzpicture}
			\begin{tikzpicture}
				[x=1.0cm,y=1.0cm]
				\clip(-0.5,-1.6) rectangle (7.5,1.6);
				\draw [color=blue] (3.6,1) -- (5.3,1.5) -- (7,1) -- (7,-1) -- (5.3,-1.5) -- (3.6,-1);
				\draw [fill=black] (0,0) circle (1pt);
				\draw [fill=black] (1.5,0) circle (1pt);
				\draw [fill=black] (3,0) circle (1pt);
				\draw [fill=black] (4.5,0) circle (1pt);
				\draw [fill=black] (6,0) circle (1pt);
				\draw [fill=black] (1.2,1) circle (1pt);
				\draw [fill=black] (2.4,1) circle (1pt);
				\draw [fill=black] (3.6,1) circle (1pt);
				\draw [fill=black] (4.8,1) circle (1pt);
				\draw [fill=black] (1.2,-1) circle (1pt);
				\draw [fill=black] (2.4,-1) circle (1pt);
				\draw [fill=black] (3.6,-1) circle (1pt);
				\draw [fill=black] (4.8,-1) circle (1pt);
				\draw [fill=black] (5.3,1.5) circle (1pt);
				\draw [fill=black] (7,1) circle (1pt);
				\draw [fill=black] (7,-1) circle (1pt);
				\draw [fill=black] (5.3,-1.5) circle (1pt);
				\draw[dash pattern=on 3pt off 3pt] (0,0) -- (1.2,1) -- (2.4,1) -- (3.6,1);
				\draw (3.6,1) -- (4.8,1) -- (6,0);
				\draw (0,0) -- (1.2,-1) -- (2.4,-1) -- (3.6,-1) -- (4.8,-1) -- (6,0);
				\draw (0,0) -- (1.5,0);
				\draw (1.5,0) -- (3,0);
				\draw (3,0) -- (4.5,0);
				\draw (4.5,0) -- (6,0);
				\node at (-0.2,0) {$s$};
				\node at (6.2,0) {$t$};
				\node at (3.6,1.2) {$s_0$};
				\node at (3.6,-1.2) {$t_0$};
				\node at (2.7,0.8) {$\sP_1$};
				\node at (2.7,-0.2) {$\sP_2$};
				\node at (2.7,-1.2) {$\sP_3$};
				\begin{scriptsize}
					\node at (7.225,0) {\textcolor{blue}{$\sP_0$}};
				\end{scriptsize}
				\node at (1.8,1.15) {$*$};
				\node at (4.2,-0.85) {$*$};
			\end{tikzpicture}
			\begin{tikzpicture}
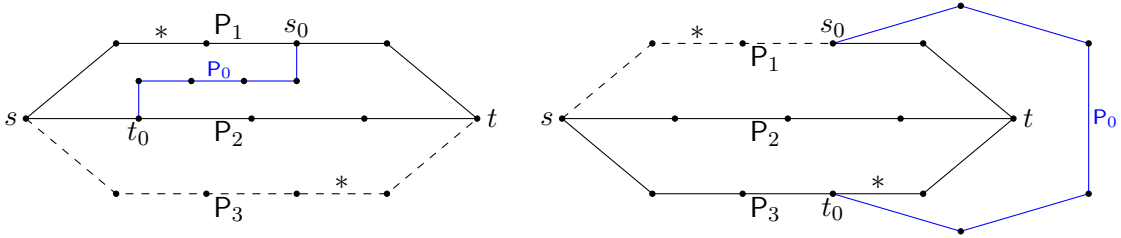
				\node at (0,0) {\textbf{\hypertarget{figfive}{Figure 5:}} Path-removal operations where $*$ indicates the repeated color. };
			\end{tikzpicture}
		\end{center}
		
		The casework above verifies our claim, and so the proof is complete. 
	\end{proof}
	
	Let $\fF$ be a Frankenstein graph with $\cP(\fF) = \{\sG_1, \dotsc, \sG_m\}$. To understand its structure better, we associate with it an auxiliary uncolored bipartite graph $G(\fF) \eqdef (V_1 \cup V_2, E)$, in which 
	\begin{itemize}
		\item $V_1 \eqdef \{\sG_1, \dotsc, \sG_m\}$, $V_2 \eqdef \{\text{the unique common vertex of some $\sG_i, \sG_j \, (i \neq j)$}\}$, and 
		\item $E \eqdef \{(\sG, v) \in V_1 \times V_2 : v \in V(\sG)\}$. 
	\end{itemize}
	
	\begin{lemma} \label{lem:Fr_aux}
		$G(\fF)$ is acyclic for every Frankenstein graph $\fF$, and so is a forest. 
	\end{lemma}
	
	\begin{proof}
		Assume to the contrary that $v_1 \sG_1 v_2 \sG_2 v_3 \dotsb v_k \sG_k v_1$ presents a cycle in $G(\fF)$, without loss of generality. Here the notations $u \sG_i$ and $\sG_j v$ refer to edges of $G(\fF)$. From \Cref{obs:bflip} we deduce that there exists for each $i \in [k]$ a rainbow path $\sP_i$ with terminals $v_i, v_{i+1}$ in $\sG_i \, (v_{k+1} = v_1)$. Since different parts in $\fF$ are edge-disjoint and color-disjoint, $\sQ \eqdef \sP_1 \cup \dotsb \cup \sP_k$ is a rainbow circuit, and so there exists a rainbow cycle $\sC \subseteq \sQ$. Since $\sC$ cannot be a subgraph of any part of $\fF$, we can find $uv \in E(\sG_x) \cap E(\sC)$ and $vw \in E(\sG_y) \cap E(\sC)$ where $x \ne y$. It follows from \ref{Fr_tree} that either $\sG_x$ or $\sG_y$, say $\sG_x$, is not a rainbow tree. However, \Cref{lem:ectest} then implies the existence of a rainbow even cycle subgraph in $\sC \cup \sG_j$, which contradicts \ref{Fr_ecfree}. 
	\end{proof}
	
	\Cref{lem:ectest,lem:Fr_aux} will be applied not only in the proof of \Cref{thm:fsteps}, but also later in many places. 
	
	\begin{proof}[Proof of \Cref{thm:fsteps}]
		We induct on $m$. The theorem is vacuously true when $m = 1$. Suppose $m \ge 2$ and let $w$ be a leaf vertex of $G(\fF)$. (If no leaf exists, then $E = \varnothing$ and any permutation $\sigma$ satisfies the theorem.) It is easily seen from the definition that no leaf exists in $V_2$, and hence we assume without loss that $w = \sG_m$. Since the partition $\{\sG_1, \dotsc, \sG_{m-1}\}$ defines a Frankenstein graph as well, the inductive hypothesis on $m-1$ implies the existence of a permutation $\sigma$ on $[m-1]$ satisfying $|V(\fF_i) \cap V(\sG_{\sigma(i+1)})| \le 1$ for all $i \in [m-2]$. Then $\sG_m$ is a leaf implies that $|V(\fF_{m-1}) \cap V(\sG_m)| \le 1$. So, by defining $\sigma(m) \eqdef m$ to extend the definition of $\sigma$, the inductive proof is complete. 
	\end{proof}
	
	The following corollaries of \Cref{thm:fsteps} will be useful in the proof of \Cref{thm:rainbow_ec}. 
	
	\begin{corollary} \label{coro:fgraph}
		If $\fF$ is a \fw graph, then $|\chi(\fF)| \leq \frac{6}{5}(|V(\fF)|-1)$. 
	\end{corollary}
	
	\begin{corollary} \label{coro:cycleinpart}
		If $\fF$ is a \fw graph with $\cP(\fF) = \{\sG_1, \dotsc, \sG_m\}$ and $\sC \subseteq \fF$ is a cycle, then there exists $i \in [m]$ such that $\sC \subseteq \sG_i$. 
	\end{corollary}
	
	\begin{proof}
		Write $V \eqdef V(\fF)$. We prove \Cref{coro:fgraph,coro:cycleinpart} by induction on $m$. 
		
		If $m=1$, then \Cref{coro:cycleinpart} is trivially true. To see that \Cref{coro:fgraph} holds, we need to check the cases when $\fF$ is a long rainbow odd cycle or a bad piece or a rainbow tree. Indeed, we have
		\[
		\begin{cases}
			|\chi(\fF)| = |V| < \frac{6}{5}(|V|-1) \quad &\text{when $\fF$ is a long rainbow odd cycle (hence $|V| \ge 7$)}, \\
			|\chi(\fF)| \le \frac{6}{5}|(|V|-1) \quad &\text{when $\fF$ is a bad piece (by \Cref{obs:bestimate})}, \\
			|\chi(\fF)| = |V|-1 < \frac{6}{5}(|V|-1) \quad &\text{when $\fF$ is a rainbow tree}. 
		\end{cases}
		\]
		
		Suppose $m \ge 2$ then. Assume without loss that the identity $\sigma(i) \eqdef i$ satisfies \Cref{thm:fsteps}. Then
		\[
		|\chi(\fF)| = |\chi(\fF_{m-1} \cup \sG_m)| = |\chi(\fF_{m-1})|+|\chi(\sG_m)| \le \frac{6}{5}(|V(\fF_{m-1})|+|V(\sG_m)|-2) \le \frac{6}{5}(|V|-1), 
		\]
		by applying the inductive hypothesis to $\fF_{m-1}$ and noticing that $|V(\fF_{m-1}) \cap V(\sG_m)| \le 1$. Also, we have $\sC \subseteq \fF_{m-1}$ or $\sC \subseteq \sG_m$ because the shared vertex of $\fF_{m-1}$ and $\sG_m$, if exists, is a cut vertex of $\fF$. By applying the inductive hypothesis to $\fF_{m-1}$, we can find some $i \in [m]$ such that $\sC \subseteq \sG_i$. 
	\end{proof}
	
	To prove \Cref{thm:rainbow_ec}, we need another technical result on Frankenstein graphs. 
	
	\begin{proposition} \label{prop:fpath}
		Suppose $\fF$ is a \fw graph and $\sP \subseteq \fF$ is a path with terminals $s$ and $t$. Then there exists a rainbow path $\sP' \subseteq \fF$ with the same terminals $s$ and $t$. 
	\end{proposition}
	
	\begin{proof}
		The existence of $\sP$ implies that $s, t$ are in the same connected component of $\fF$. We thus assume without loss that $\fF$ is connected. Then there exists a path in the uncolored graph $G(\fF)$ of the form $\sG_{i_1} v_1 \sG_{i_2} v_2 \dotsb v_{\ell-1} \sG_{i_{\ell}}$
		such that $s \in V(\sG_{i_1}), \, t \in V(\sG_{i, \ell})$ and $\ell \ge 1$. It then follows from \Cref{obs:bflip} that there exists a rainbow trail $\sQ$ joining $s$ and $t$. Obviously, any path $\sP' \subseteq \sQ$ with terminals $s$ and $t$ satisfies \Cref{prop:fpath}. 
	\end{proof}
	
	For a Frankenstein graph $\fF$ given by the partition $\cP(\fF) = \{\sC_1, \dotsc, \sC_c, \sB_1, \dotsc, \sB_b, \sT_1, \dotsc, \sT_t\}$, we associate with it counting parameters $c(\fF) \eqdef c$ and $b(\fF) \eqdef b$. Notice that $c(\fF), b(\fF)$ depend on not only the graph $\fF$, but the partition $\cP(\fF)$ as well. We still need another depth parameter. 
	
	For any tree $\sT$ with $V(\sT) \subset \N_+$, let its \emph{root} be $r \eqdef \min V(\sT)$. For any vertex $v \in V(\sT)$, define its \emph{relative depth} in $\sT$ as $\depth_{\sT}(v) \eqdef \dist_{\sT}(r, v)$, which is the length of the unique path with terminals $r$ and $v$. We henceforth define for any forest $\sF$ with $V(\sF) \subset \N_+$ its \emph{total depth} as 
	\[
	\Depth(\sF) \eqdef \sum_{i=1}^t \sum_{v \in V(\sT_i)} \depth_{\sT_i}(v) 
	\]
	where $\sT_1, \dotsc, \sT_t$ are the connected components of $\sF$. For any \fw graph $\fF$ with $V(\fF) \subset \N_+$, we refer to its \emph{total depth} as the total depth of its forest part, i.e.~$\Depth(\fF) \eqdef \Depth(\sT_1 \cup \dotsb \cup \sT_t)$. 
	
	Later in practice, we shall often construct a Frankenstein graph by a ``partition'' 
	\[
	\cP(\fF) = \{\sC_1, \dotsc, \sC_c, \sB_1, \dotsc, \sB_b, \sF\}
	\]
	where $\sC$'s are long rainbow odd cycles, $\sB$'s are bad pieces, $\sF = \sT_1 \cup \dotsb \cup \sT_t$ is the union of vertex-disjoint and color-disjoint rainbow trees, such that $\chi(\sG_i) \cap \chi(\sG_j) = \varnothing$ for any distinct $\sG_i, \sG_j \in \cP(\fF)$. Indeed, this $\cP(\fF)$ is formally not a partition since $\sF$ and $\sC_i$ or $\sB_j$ may share more than one vertex. However, \ref{Fr_tree} implies, up to a relabeling of the rainbow tree parts of $\fF$, that there is no difference between exposing the trees $\sT_1, \dotsc, \sT_t$ and exposing the forest $\sF$. 
	
	\section{Proof of \texorpdfstring{\Cref{thm:rainbow_ec}}{Theorem~\ref{thm:rainbow_ec}}} \label{sec:proof_ec}
	
	We prove \Cref{thm:rainbow_ec} indirectly. Suppose $\cD = (\sD_1, \dotsc, \sD_m)$ is a family of $m \eqdef \lf\frac{6(n-1)}{5}\rf+1 > \frac{6(n-1)}{5}$ even cycles on the ambient vertex set $[n]$ without any rainbow even cycle subgraph. 
	
	\smallskip
	
	Let $\fF_*$ be a \fw subgraph of the family $\cD$ satisfying the following maximal conditions: 
	\begin{enumerate}[label=(M\arabic*),  ref=(M\arabic*)]
		\item\label{max:cycle} The number of long rainbow odd cycles $c(\fF_*)$ is maximized. 
		\item\label{max:bpiece} The number of bad pieces $b(\fF_*)$ is maximized under \ref{max:cycle}. 
		\item\label{max:edges} The number of edges $|\fF_*|$ is maximized under \ref{max:bpiece}. 
		\item\label{min:depth} The total depth $\Depth(\fF_*)$ is minimized under \ref{max:edges}. 
	\end{enumerate}
	
	Suppose the partition of $\fF_*$ is 
	\[
	\cP(\fF_*) = \{\sC_1, \dotsc, \sC_c, \sB_1, \dotsc, \sB_b, \sT_1, \dotsc, \sT_t\} \quad \text{with} \quad \sF \eqdef \sT_1 \cup \dotsb \cup \sT_t, 
	\]
	where $\sC$'s are long rainbow odd cycles, $\sB$'s are bad pieces, and $\sT$'s are vertex-disjoint rainbow trees. 
	
	\subsection{Outer edges and outer cycles}
	
	Let $\lambda$ be the color of the even cycle $\sD_{\lambda}$. From \Cref{coro:fgraph} we deduce that $|\chi(\fF_*)| \le \frac{6}{5}(n-1) < |\cD|$, and hence $\Lambda \eqdef [m] \setminus \chi(\fF_*) \ne \varnothing$. Indeed, every edge of the multigraph $\cD_{\Lambda} \eqdef \bigcup_{\lambda \in \Lambda} \sD_{\lambda}$ is absent in $\fF_*$. 
	
	We call $\sf$ in $\cD_{\Lambda}$ an \emph{outer edge} if no coincident edge of $\sf$ is in $\fF_*$. A rainbow $\{3, 5\}$-cycle containing an outer edge $\sf$ in $\fF_* + \sf$ is called an \emph{outer cycle} of $\sf$. Hereafter $\sG + \se$ denotes the graph generated by adding $\se$ to $\sG$ (i.e.~$\sG + \se \eqdef \sG \cup \{\se\}$). Moreover, whenever we write $\sG + \se$, we implicitly assume that $\se$ is not coincident to any edge of $\sG$. Similarly, $\sG - \se$ (assuming $\se \in \sG$) refers to the graph obtained by deleting $\se$ from $\sG$ (i.e.~$\sG - \se \eqdef \sG \setminus \{\se\}$). Recall that a (colored) graph is a set of colored edges. 
	
	\smallskip
	
	The next propositions are devoted to the existence of outer edges and outer cycles. 
	
	\begin{proposition} \label{prop:outeredge}
		For any $\lambda \in \Lambda$, an outer edge exists in $\sD_{\lambda}$. 
	\end{proposition}
	
	\begin{proof}
		Assume for the sake of contradiction that $\sD_{\lambda}$ is covered by $\fF_*$. That is, each $\se \in \sD_{\lambda}$ has one coincident edge $\se^* \in \fF_*$. Indeed, this $\se^*$ is unique because no coincident edges exist in a Frankenstein graph. Define $\sD_{\lambda}^* \eqdef \{\se^* : \se \in \sD_{\lambda}\} \subseteq \fF_*$. Since long rainbow odd cycles and rainbow trees contain no even cycle, it follows from \Cref{coro:cycleinpart} that $\sD_{\lambda}^*$ has to be contained in some bad piece $\sB_j$, and so $|\sD_{\lambda}^*| - |\chi(\sD_{\lambda}^*)| \in \{0, 1\}$. Since no rainbow even cycle exists in $\cD$, we obtain $|\sD_{\lambda}^*| - |\chi(\sD_{\lambda}^*)| = 1$. So, there exists a unique pair of distinct edges $\se_1^*, \se_2^*$ in $\sD_{\lambda}^*$ such that $\chi(\se_1^*) = \chi(\se_2^*)$. Thus, $\sD_{\lambda}^* - \se_1^* + \se_1$ is a rainbow even cycle in $\cD$, a contradiction. 
	\end{proof}
	
	\begin{proposition} \label{prop:outercycle}
		For any outer edge $\sf$, an outer cycle of $\sf$ exists. 
	\end{proposition}
	
	\begin{proof}
		Let $V(\sf) \eqdef \{u, v\}$. Observe that $u, v$ are in a same connected component of $\fF_*$, for otherwise
		\[
		\cP(\fF_* + \sf) \eqdef \{ \sC_1, \dotsc, \sC_c, \sB_1, \dotsc, \sB_b, \sF+\sf \} 
		\]
		gives another \fw subgraph of $\cD$ with one more edge than $\fF_*$, which contradicts \ref{max:edges}. It follows from \Cref{prop:fpath} that $\sf$ completes a rainbow (hence odd) cycle $\sC^{\sf}$ in $\sF+\sf$. 
		
		It then suffices to disprove that $\sC^{\sf}$ is long. Assume to the contrary that $\sC^{\sf}$ is long. Since $\sf \notin \fF_*$, from \Cref{lem:ectest} we deduce that $|V(\sC^{\sf}) \cap V(\sC_i)| \le 1$ for every $i \in [c]$. So, $\mathcal{P}(\fF_{+}) \eqdef \{\sC_1, \dotsc, \sC_c, \sC^{\sf}\}$ presents another Frankenstein subgraph of $\cD$ with $c(\fF_{+}) > c(\fF_*)$, which contradicts \ref{max:cycle}. 
	\end{proof}
	
	For any tree $\sT$ with $v \in V(\sT) \subseteq [n]$, we define 
	\[
	\child_{\sT}(v) \eqdef \bigl\{ w \in V(\sT) : vw \in E(\sT), \, \depth_{\sT}(w) = \depth_{\sT}(v) + 1 \bigr\}. 
	\]
	The following properties characterize behaviors of outer $3$-cycles. 
	
	\begin{proposition} \label{prop:outerthree}
		Suppose $\sf$ is an outer edge with $V(\sf) = \{u, v\}$, and $\sC$ is an outer $3$-cycle of $\sf$ with $V(\sC) = \{u, v, w\}$. Then there exists $k \in [t]$ such that $u, v, w \in V(\sT_k)$ and $u, v \in \child_{\sT_k}(w)$. 
	\end{proposition}
	
	\begin{proof}
		We show $u, v, w \in V(\sT_k)$ for some $k$ first. Since $uw, vw \in E(\fF_*)$, we may assume $uw \in E(\sX_1)$ and $vw \in E(\sX_2)$ where $\sX_1, \sX_2 \in \cP(\fF_*)$. In fact, $\sX_{\bullet} \in \{\sT_1, \dotsc, \sT_t\} \, (\bullet = 1, 2)$, for otherwise \Cref{lem:ectest} implies the existence of a rainbow even cycle in $\sC \cup \sX_{\bullet}$. It follows from \ref{Fr_tree} that $\sX_1 = \sX_2 = \sT_k$. 
		
		We prove $u, v \in \child_{\sT_k}(w)$ then. Suppose $\se_1 \eqdef (uw, \alpha)$ and $\se_2 \eqdef (vw, \beta)$ are edges in $\sT_k$. The existence of $\se_1, \se_2$ tells us that $\abs{\depth_{\sT_k}(u) - \depth_{\sT_k}(v)}$ is either $0$ or $2$. It suffices to establish that $\depth_{\sT_k}(u) = \depth_{\sT_k}(v)$. If not, then assume without loss that $\depth_{\sT_k}(u) = \depth_{\sT_k}(v) + 2$. Since $\depth_{\sT_k'}(u) < \depth_{\sT_k}(u)$ and $\depth_{\sT'_k}(x) \le \depth_{\sT_k}(x)$ for all $x \in V(\sT_k) = V(\sT'_k)$, we deduce that $\sT'_k \eqdef \sT_k+\sf-\se_1$ is another tree with $\Depth(\sT'_k) < \Depth(\sT_k)$. Then the partition 
		\[
		\cP(\fF') \eqdef \cP(\fF_*+\sf-\se_2) = \{ \sC_1, \dotsc, \sC_c, \sB_1, \dotsc, \sB_b, \sT_1, \dotsc, \sT_k', \dotsc, \sT_t \}
		\]
		gives a \fw subgraph of $\cD$. However, this contradicts \ref{min:depth} since $\Depth(\fF') < \Depth(\fF_*)$. Therefore, $\depth_{\sT_k}(u) = \depth_{\sT_k}(v)$, and so $u, v \in \child_{\sT_k}(w)$. 
	\end{proof}
	
	\begin{proposition} \label{prop:onlyouterthree}
		Suppose no outer $5$-cycle exists. If $\sf = (uv, \alpha)$ is an outer edge with outer cycle $\sC$ on vertices $u, v, w \in \sT_k$ (by \Cref{prop:outerthree}), then $\sD_{\alpha}$, the even cycle of color $\alpha$ from $\cD$ containing $\sf$, satisfies $V(\sD_{\alpha}) \subseteq \{w\} \cup \child_{\sT_k}(w)$. (See \hyperlink{figsix}{Figure~6}.)
	\end{proposition}
	
	\begin{center}
		\begin{tikzpicture}
			\draw [fill=black] (-2,1) circle (1pt);
			\draw [fill=black] (-4,0) circle (1pt);
			\draw [fill=black] (0,0) circle (1pt);
			\draw [fill=black] (-5,-1) circle (1pt);
			\draw [fill=black] (-4,-1) circle (1pt);
			\draw [fill=black] (-3,-1) circle (1pt);
			\draw [fill=black] (-1.5,-1) circle (1pt);
			\draw [fill=black] (-0.5,-1) circle (1pt);
			\draw [fill=black] (0.5,-1) circle (1pt);
			\draw [fill=black] (1.5,-1) circle (1pt);
			\draw [fill=black] (-5.4,-2) circle (1pt);
			\draw [fill=black] (-5,-2) circle (1pt);
			\draw [fill=black] (-4.6,-2) circle (1pt);
			\draw [fill=black] (-3.2,-2) circle (1pt);
			\draw [fill=black] (-2.8,-2) circle (1pt);
			\draw [fill=black] (-1.7,-2) circle (1pt);
			\draw [fill=black] (-1.3,-2) circle (1pt);
			\draw [fill=black] (0.1,-2) circle (1pt);
			\draw [fill=black] (0.5,-2) circle (1pt);
			\draw [fill=black] (0.9,-2) circle (1pt);
			\draw [fill=black] (1.5,-2) circle (1pt);
			\draw (-2,1) -- (-4,0);
			\draw (-2,1) -- (0,0);
			\draw (-4,0) -- (-5,-1);
			\draw (-4,0) -- (-4,-1);
			\draw (-4,0) -- (-3,-1);
			\draw (0,0) -- (-1.5,-1);
			\draw (0,0) -- (-0.5,-1);
			\draw (0,0) -- (0.5,-1);
			\draw (0,0) -- (1.5,-1);
			\draw (-5,-1) -- (-5.4,-2);
			\draw (-5,-1) -- (-5,-2);
			\draw (-5,-1) -- (-4.6,-2);
			\draw (-3,-1) -- (-3.2,-2);
			\draw (-3,-1) -- (-2.8,-2);
			\draw (-1.5,-1) -- (-1.7,-2);
			\draw (-1.5,-1) -- (-1.3,-2);
			\draw (0.5,-1) -- (0.1,-2);
			\draw (0.5,-1) -- (0.5,-2);
			\draw (0.5,-1) -- (0.9,-2);
			\draw (1.5,-1) -- (1.5,-2);
			\draw[color=blue] (-0.5,-1) -- (0.5,-1);
			\draw[color=black] (-0.7,-1) node {$u$};
			\draw[color=black] (0.7,-1) node {$v$};
			\draw[color=black] (0.2,0.1) node {$w$};
			\draw[color=blue] (0,-0.8) node {$\sf$};
			\draw[dash pattern=on 3pt off 3pt] (-1.7,-1.25) rectangle (1.7,0.25);
			\draw[dash pattern=on 3pt off 3pt] (-5.9,-1.5) rectangle (1.9,0.5);
			\draw[dash pattern=on 3pt off 3pt] (-5.9,1.2) -- (-5.9,0.5);
			\draw[dash pattern=on 3pt off 3pt] (1.9,1.2) -- (1.9,0.5);
			\draw[dash pattern=on 3pt off 3pt] (-5.9,-2.2) -- (-5.9,-1.5);
			\draw[dash pattern=on 3pt off 3pt] (1.9,-2.2) -- (1.9,-1.5);
			\node at (-5.6,0.75) {$A^-$};
			\node at (-5.65,0.25) {$A'$};
			\node at (-5.6,-1.7) {$A^+$};
			\node at (-1.5,0.05) {$A$};
			\node at (-2, -2.7) {\textbf{\hypertarget{figsix}{Figure 6:}} $V(\sf) \subseteq \child_{\sT_k}(w)$ implies $V(\sD_{\alpha}) \subseteq A$. };
		\end{tikzpicture}
	\end{center}
	
	\vspace{-2em}
	\begin{proof}
		We first show that $V(\sD_{\alpha}) \subseteq V(\sT_k)$. Define $\tau \colon \sD_{\alpha} \to \cP(\fF_*)$ as follows: For any edge $\se \in \sD_{\alpha}$, 
		\vspace{-0.5em}
		\begin{itemize}
			\item if $\se$ is an outer edge, then $V(\se) \subseteq V(\sT_{\ell})$ for some $\ell$ (by \Cref{prop:outerthree}), and we set $\tau(\se) \eqdef \sT_{\ell}$; 
			\vspace{-0.5em}
			\item if $\se$ is coincident to $\se' \in \fF_{*}$, then we set $\tau(\se) \eqdef \sX$, where $\sX$ is the part of $\fF_{*}$ that contains $\se'$. 
		\end{itemize}
		\vspace{-0.5em}
		By applying $\tau$ on $\sD_{\alpha}$, we locate a closed walk $Q \subseteq G(\fF_*)$ as follows: 
		\vspace{-0.5em}
		\begin{enumerate}
			\item[(1)] Put the edges of $\sD_{\alpha}$ on a circle $\mathcal{O}$ in order. Replace $\se$ by $\tau(\se)$ for each $\se \in \sD_{\alpha}$. 
			\vspace{-0.5em}
			\item[(2)] If two consecutive objects on $\mathcal{O}$ are the same, then remove one of them. Repeat. 
			\item[(3)] If $\sG_i, \sG_j \in \cP(\fF_*)$ are adjacent on $\mathcal{O}$, then plug in $v_{ij} \in V(\sG_i) \cap V(\sG_j)$ between them. 
		\end{enumerate}
		\vspace{-0.25em}
		The resulting arrangement on $\mathcal{O}$ forms a closed walk $Q \subseteq G(\fF_*)$, where each pair of consecutive edges $v_{ij} \sG_j, \, \sG_j v_{jk}$ in $Q$ correspond to a path with terminals $v_{ij}, v_{jk}$ on $\sD_{\alpha}$. Indeed, $Q$ is a circuit because $\sD_{\alpha}$ passes through each $v_{ij}$ exactly once. For instance, if $\sD_{\alpha}$ consists of $\se_1, \dotsc, \se_8$ in order such that 
		\[
		\tau(\se_1, \se_2, \se_3, \se_4, \se_5, \se_6, \se_7, \se_8) = (\sG_3, \sG_3, \sG_2, \sG_2, \sG_1, \sG_3, \sG_4, \sG_5), 
		\]
		then the steps (1) through (3) generate 
		
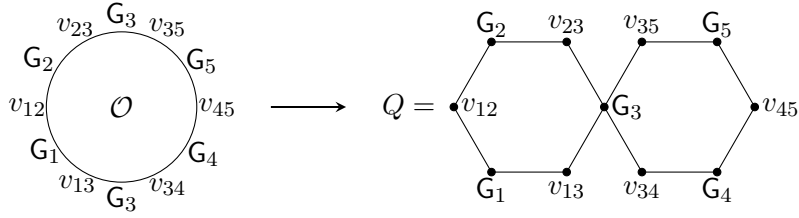
\[
		\begin{tikzpicture}[scale = 0.5]
			\clip (-3,-3) rectangle (3,3);
			\draw (0,0) circle (2);
			\node at (0,0) {$\mathcal{O}$};
			\node at (90:2.45) {$\sG_3$};
			\node at (120:2.45) {$v_{23}$};
			\node at (150:2.55) {$\sG_2$};
			\node at (180:2.5) {$v_{12}$};
			\node at (210:2.35) {$\sG_1$};
			\node at (240:2.35) {$v_{13}$};
			\node at (270:2.45) {$\sG_3$};
			\node at (300:2.45) {$v_{34}$};
			\node at (330:2.5) {$\sG_4$};
			\node at (0:2.55) {$v_{45}$};
			\node at (30:2.45) {$\sG_5$};
			\node at (60:2.45) {$v_{35}$};
		\end{tikzpicture}
		\quad
		\begin{tikzpicture}[scale = 0.5]
			\clip(-1.2,-3) rectangle (1.1,3);
			\draw[line width = 0.75pt][-stealth] (-1,0) -- (1,0);
		\end{tikzpicture}
		\quad 
		\begin{tikzpicture}[scale = 0.5]
			\clip (-6,-3) rectangle (5.2,3);
			\draw (0,0) -- (-1,1.732) -- (-3,1.732) -- (-4,0) -- (-3,-1.732) -- (-1,-1.732) -- (0,0) -- (1,1.732) -- (3,1.732) -- (4,0) -- (3,-1.732) -- (1,-1.732) -- (0,0);
			\foreach \x in {-3,-1,1,3}
			\foreach \y in {-1.732,1.732}
			\draw[fill=black] (\x,\y) circle (0.1);
			\foreach \x in {-4,0,4}
			\draw[fill=black] (\x,0) circle (0.1);
			\node at (0.6,0) {$\sG_3$};
			\node at (-1,-2.2) {$v_{13}$};
			\node at (-1,2.2) {$v_{23}$};
			\node at (1,-2.2) {$v_{34}$};
			\node at (1,2.2) {$v_{35}$};
			\node at (-3.3,0) {$v_{12}$};
			\node at (4.7,0) {$v_{45}$};
			\node at (-3,2.2) {$\sG_2$};
			\node at (-3,-2.2) {$\sG_1$};
			\node at (3,2.2) {$\sG_5$};
			\node at (3,-2.2) {$\sG_4$};
			\node at (-5.2,0) {$Q =$};
		\end{tikzpicture}
		\]
		
		\vspace{-0.75em}
		\noindent However, \Cref{lem:Fr_aux} asserts that $\fF_*$ is acyclic. So, $Q$ is a single vertex, and hence $V(\sD_{\alpha}) \subseteq V(\sT_k)$. 
		
		Use abbreviations $V \eqdef V(\sT_k)$ and $d(x) \eqdef \depth_{\sT_k}(x)$. Partition $V$ into $A \eqdef \{w\} \cup \child_{\sT_k}(w)$, $A^+ \eqdef \{x \in V: d(x) > d(v)\}$, $A^- \eqdef \{x \in V: d(x) < d(w)\}$ and $A' \eqdef V \setminus (A \cup A^+ \cup A^-)$. Let $T_k$ be the uncolored copy of $\sT_k$, which is the uncolored graph on vertex set $V = V(\sT_k)$ and edge set $E(\sT_k)$. For all $z \in V$ and all pairs of distinct vertices $x, y \in \child_{T_k}(z)$, we add the edges $xy$ simultaneously into $E(T_k)$ to generate a new graph $\overline{T}_k$. The vertex set of $\overline{T}_k$ is still $V$. Due to the absence of outer $5$-cycles, from \Cref{prop:outerthree} we deduce that $D_{\alpha}$, the uncolored copy of $\sD_{\alpha}$, is a subgraph of $\overline{T}_k$. 
		
		Notice that any subpath of $T_k$ with one terminal in $A$ and the other in $A'$ must go through $A^-$. It then suffices to show that $V(D_{\alpha}) = V(\sD_{\alpha})$ and $A^-$ are disjoint. This breaks down to exclude the situation $d(z_+) \ge d(z_-)+2$ for some $z_+, z_- \in V(D_{\alpha})$. If such $z_+, z_-$ exist, then $D_{\alpha}$ consists of two subpaths $P_1, P_2$ with terminals $z_+$ and $z_-$. Let $z_i$ be the vertex on $P_i$ with $d(z_i) = d(z_+)-1$ that is nearest to $z_+$. The crucial observation is that $z_i$ is the parent of $z_+$, which is the unique vertex in $V$ such that $z_+ \in \child_{T_k}(z_i)$. Indeed, this follows from the fact that $z_i$ is a cut vertex of $\overline{T}_k$ which separates $z_+$ from all vertices of smaller depths. However, the observation implies that $z_1 = z_2$, which is absurd. We conclude that $V(\sD_{\alpha}) = V(D_{\alpha}) \subseteq A$, and hence the proof is complete. 
	\end{proof}
	
	\subsection{Finishing the proof}
	
	\begin{lemma} \label{lem:f'_5cycle}
		There exists a \fw subgraph $\fF_0$ of $\cD$ whose partition is given by 
		\[
		\cP(\fF_0) = \{\sC_1, \dotsc, \sC_c, \sB_1, \dotsc, \sB_b, \sF_0\} \quad \text{with} \quad |\sF_0| = |\sF|, 
		\]
		and an edge $\sf_0$ in $\cD$ such that
		$\chi(\sf_0) \notin \chi(\fF_0)$ and $\sf_0$ completes a rainbow $5$-cycle in $\fF_0+\sf_0$. 
	\end{lemma}
	
	\begin{proof}
		If there is an outer $5$-cycle in $\fF_*$, say $\sC^{\overline{\sf}}$ of an outer edge $\overline{\sf}$, then $(\fF_0, \sf_0) \eqdef (\fF_*, \overline{\sf})$ with $\sF_0 \eqdef \sF$ satisfies \Cref{lem:f'_5cycle}. We assume no outer $5$-cycle exists then. It follows from \Cref{prop:outeredge,prop:outercycle} that an outer edge $\sf$ and its outer cycle $\sC^{\sf}$ exist. Suppose $\sf \eqdef (uv, \alpha)$ and $\sD_{\alpha}$ is the monochromatic even cycle from $\cD$ that contains $\sf$. Assume $V(\sC^{\sf}) \eqdef \{u, v, w\}$. It follows from \Cref{prop:outerthree} that $u, v, w$ all lie in a single rainbow tree $\sT_k \in \cP(\fF_*)$ and $u, v \in \child_{\sT_k}(w)$. From \Cref{prop:onlyouterthree} we deduce that $V(\sD_{\alpha}) \subseteq \{w\} \cup \child_{\sT_k}(w)$. Since $\sD_{\alpha}$ consists of at least $4$ edges, at least $1$ of the two adjacent edges of $\sf$ on $\sD_{\alpha}$ is not incident to the vertex $w$. Assume without loss that $\sf' \eqdef (uv', \alpha)$ is such an edge, and hence $v' \in \child_{\sT_k}(w)$. Observe that $v$ and $v'$ are symmetric despite our definition. 
		
		\begin{center}
			\begin{tikzpicture}[scale = 1.8]
				\draw (-2,0.025) -- (0,1.025) -- (2,0.025);
				\draw (-1,0) -- (0,1);
				\draw (1,0) -- (0,1);
				\draw[color=green] (0,0) -- (-1,0) -- (-2,0) -- (0,1) -- (2,0) -- (1,0) -- (0,0);
				\draw[color=blue] (0,0) -- (0,1) -- (-1,1.5) -- (-2,1.5) -- (-3,1) -- (-3,0) -- (-2,-0.5) -- (-1,-0.5) -- (0,0);
				\draw[color=red] (-3,1) -- (0,1);
				\draw[color=red] (-3,0) arc[start angle=210,end angle=330,radius=1.732cm];
				\draw[fill=black] (-2,0) circle (1.5pt);
				\draw[fill=black] (-1,0) circle (1.5pt);
				\draw[fill=black] (0,0) circle (1.5pt);
				\draw[fill=black] (1,0) circle (1.5pt);
				\draw[fill=black] (2,0) circle (1.5pt);
				\draw[fill=black] (0,1) circle (1.5pt);
				\draw[fill=black] (-1,1.5) circle (1.5pt);
				\draw[fill=black] (-2,1.5) circle (1.5pt);
				\draw[fill=black] (-3,1) circle (1.5pt);
				\draw[fill=black] (-3,0) circle (1.5pt);
				\draw[fill=black] (-2,-0.5) circle (1.5pt);
				\draw[fill=black] (-1,-0.5) circle (1.5pt);
				\node at (0.27,-0.15) {$x_0 = u$};
				\node at (0.45,1.1) {$x_{2k+1} = w$};
				\node at (-1,-0.125) {$v'$};
				\node at (1,-0.15) {$v$};
				\node at (-1,-0.65) {$x_1$};
				\node at (-1,1.65) {$x_{2k}$};
				\node at (-2,-0.65) {$\cdots$};
				\node at (-2,1.65) {$\cdots$};
				\node at (-3.2,-0.125) {$x_{\lambda-1}$};
				\node at (-3.15,1.125) {$x_{\lambda}$};
				\node at (0.5,0.1) {\textcolor{green}{$\sf$}};
				\node at (-0.47,0.11) {\textcolor{green}{$\sf'$}};
				\node at (1,0.4) {\textcolor{green}{$\sD_{\alpha}$}};
				\node at (-0.075,0.5) {\textcolor{blue}{$\sg$}};
				\node at (0.4,0.5) {$\sk$};
				\node at (-2.925,0.5) {\textcolor{blue}{$\sh$}};
				\node at (-1.5,1.38) {\textcolor{blue}{$\sD_{\beta}$}};
				\node at (-1.5,0.9) {$\textcolor{red}{\sP_w}$};
				\node at (-1.5,-0.75) {$\textcolor{red}{\sP_u}$};
				\draw[color=violet,dash pattern=on 3pt off 3pt] (-3.6,0.35) -- (-2.4,0.35) -- (-2.4,-0.25) -- (-0.6,-0.25) -- (-0.6,0.35) -- (0.6,0.35) -- (0.6,-0.9);
				\node at (-3.4,0.55) {\textcolor{violet}{$\sS_w$}};
				\node at (-3.4,0.15) {\textcolor{violet}{$\sS_u$}};
				\node at (-0.8,-1.2) {\textbf{\hypertarget{figseven}{Figure 7:}} An illustration of the proof of \Cref{lem:f'_5cycle}. };
			\end{tikzpicture}
		\end{center}
		
		Let $\sg \eqdef (uw, \beta) \in \fF_*$ be the edge with $V(\sg) = \{u, w\}$. Suppose 
		\[
		\sD_{\beta} = \sg + (x_0x_1, \beta) + (x_1x_2, \beta) + \dotsb + (x_{2k}x_{2k+1}, \beta) \qquad (x_0 \eqdef u, \, x_{2k+1} \eqdef w, \, k \in \N_+)
		\]
		is the monochromatic even cycle from $\cD$ containing $\sg$. From \Cref{lem:Fr_aux} we deduce that there are two connected components $\sS_u$ and $\sS_w$ in the graph $\fF_* - \sg$ such that $u \in V(\sS_u)$ and $w \in V(\sS_w)$. Define $\lambda$ as the smallest index such that $x_{\lambda} \notin V(\sS_u)$ and write $\sh \eqdef (x_{\lambda-1}x_{\lambda}, \beta)$. Then $\sh \notin \fF_*$. 
		
		We claim that $x_{\lambda} \in \sS_w$. If not, then $\sh$ cannot complete any cycle in $\widehat{\sF} \eqdef \sF + \sf - \sg + \sh$, and so $\widehat{\sF}$ is a rainbow forest. Observe that the trees in $\widehat{\sF}$ containing $\sf$ or $\sh$ (they are possibly the same) share at most one vertex with any of $\sC_1, \dotsc, \sC_c, \sB_1, \dotsc, \sB_b$. This implies that the partition 
		\[
		\cP(\fF_* + \sf - \sg + \sh) \eqdef \{\sC_1, \dotsc, \sC_c, \sB_1, \dotsc, \sB_b, \widehat{\sF}\} 
		\]
		presents another Frankenstein subgraph of $\cD$ on $|\fF_*|+1$ edges, which contradicts \ref{max:edges}. 
		
		By \Cref{prop:fpath}, we can find a rainbow path $\sP_u \subseteq \sS_u$ with terminals $u, x_{\lambda-1}$ and a rainbow path $\sP_{w} \subseteq \sS_w$ with terminals $w, x_{\lambda}$. Here we allow $\sP_u$ to be empty if $x_0 = x_{\lambda-1}$, and allow $\sP_w$ to be empty if $x_{2k+1} = x_{\lambda}$. Note that $\sP_u = \sP_w = \varnothing$ cannot happen, since $|\sD_{\beta}| \ge 4$. Assume further that the length of $\sP_w$ is minimized, and so $v \notin V(\sP_w)$ or $v' \notin V(\sP_w)$, say $v \notin V(\sP_w)$. Thus, $\widetilde{\sC} \eqdef \sf + \sP_u +\sh + \sP_w + \sk$ is a rainbow odd cycle with $|\widetilde{\sC}| \ge 5$. Here $\sk$ denotes the edge of $\sT_k$ with $V(\sk) = \{v, w\}$. 
		
		Since $V(\sT_k+\sf-\sg) = V(\sT_k)$, we can define another Frankenstein subgraph $\fF_0 \eqdef \fF_*+\sf-\sg$ by 
		\[
		\cP(\fF_0) \eqdef \{\sC_1, \dotsc, \sC_c, \sB_1, \dotsc, \sB_b, \sT_1, \dotsc, \sT_k+\sf-\sg, \dotsc, \sT_t\}. 
		\]
		We claim that $\sf_0 \eqdef \sh$ is as desired. It suffices to show that $\widetilde{\sC}$ is a rainbow $5$-cycle. Since $\sh \in \widetilde{\sC}$ and $\beta \notin \chi(\fF_0)$, \Cref{lem:ectest} tells us that $|V(\widetilde{\sC}) \cap V(\sC_i)| \le 1 \, (\forall i \in [c])$. Then $\mathcal{P}(\fF') \eqdef \{\sC_1, \dotsc, \sC_c, \widetilde{\sC}\}$ gives another Frankenstein subgraph of $\cD$ with $c(\fF') > c$ if $\widetilde{\sC}$ is long, which contradicts \ref{max:cycle}. Thus, $|\widetilde{\sC}| \ge 5$ implies that $\widetilde{\sC}$ is a rainbow $5$-cycle in $\fF_0+\sh$. The proof of \Cref{lem:f'_5cycle} is complete. 
	\end{proof}
	
	\begin{lemma} \label{lem:5cyclegrowth}
		Suppose the rainbow $5$-cycle found in \Cref{lem:f'_5cycle} is $\widetilde{\sC} \eqdef \{(v_iv_{i+1}, \alpha_i) : i \in [5]\}$, with the convention $v_{\ell+5} = v_{\ell}$. Write $\se_i \eqdef (v_iv_{i+1}, \alpha_i)$. Then there exists a shifting parameter $j \in \{0, 1, 2, 3, 4\}$, a set of five edges $\se_i' \eqdef (v_iv_{i+1}, \alpha_{i+j})$ from $\cD$, a vertex $v^* \in [n] \setminus \{v_1, \dotsc, v_5\}$, and an index $k \in [5]$, such that at least one of the edges $(v^*v_k, \alpha_{k+j-1})$ and $(v^*v_k, \alpha_{k+j})$ appears in $\cD$. 
	\end{lemma}
	
	Informally speaking, \Cref{lem:5cyclegrowth} is dedicated to ``grow'' one more edge from the $5$-cycle guaranteed by \Cref{lem:f'_5cycle}. That is, after a possible cyclic shift of the colors on $\widetilde{\sC}$, we would like to find out another edge on one of the monochromatic even cycles in $\cD$ ``leaving'' $\widetilde{\sC}$ (i.e.~incident to $v^* \notin \{v_1, \dotsc, v_5\}$). Such a configuration will then help us to locate another bad piece in $\cD$, which contradicts \ref{max:bpiece}. For ease of notations, we write $(a, b, c, d, e) \eqdef (\alpha_1, \alpha_2, \alpha_3, \alpha_4, \alpha_5)$ in the coming example and in the proof of \Cref{lem:5cyclegrowth}. \hyperlink{figeight}{Figure~8} illustrates one possible output of \Cref{lem:5cyclegrowth} in which $(j, k) = (4, 3)$. 
	
	\begin{center}
		\begin{tikzpicture}[x=1.0cm,y=1.0cm]
			\clip(-2.5,-1.5) rectangle (2.5,2.2);
			\draw  (0.,2.)-- (-2.,1.);
			\draw  (-2.,1.)-- (-1.,-1.);
			\draw  (-1.,-1.)-- (1.,-1.);
			\draw  (1.,-1.)-- (2.,1.);
			\draw  (2.,1.)-- (0.,2.);
			\draw [fill=black] (-1.,-1.) circle (3pt);
			\draw[color=black] (-0.8,-0.8) node {$v_3$};
			\draw [fill=black] (1.,-1.) circle (3pt);
			\draw[color=black] (0.8,-0.8) node {$v_4$};
			\draw [fill=black] (-2.,1.) circle (3pt);
			\draw[color=black] (-1.7,0.9) node {$v_2$};
			\draw [fill=black] (2.,1.) circle (3pt);
			\draw[color=black] (1.7,0.9) node {$v_5$};
			\draw [fill=black] (0.,2.) circle (3pt);
			\draw[color=black] (0,1.75) node {$v_1$};
			\draw[color=black] (-0.9,1.35) node {$a$};
			\draw[color=black] (-1.325,0.075) node {$b$};
			\draw[color=black] (0,-0.8) node {$c$};
			\draw[color=black] (1.325,0.075) node {$d$};
			\draw[color=black] (0.9, 1.35) node {$e$};
		\end{tikzpicture}
		\begin{tikzpicture}
			\clip(-1.2,-1.2) rectangle (1.2,2.2);
			\draw[line width = 0.75pt][-stealth] (-1,0.5) -- (1,0.5);
			\node at (-0.05,0.9) {\Cref{lem:5cyclegrowth}};
		\end{tikzpicture}
		\begin{tikzpicture}[x=1.0cm,y=1.0cm]
			\clip(-3.2,-1.5) rectangle (2.5,2.2);
			\draw (0,1.95) -- (-1.95,1) -- (-0.95,-0.95) -- (-0.95,-0.97) -- (0.95,-0.97) -- (0.95,-0.95) -- (1.95,1) -- (0,1.95);
			\draw (0,2.05) -- (-1.95,1.1) -- (-2.05,1) -- (-1.05,-0.95) -- (-1.05,-1.05) -- (1.05,-1.05) -- (1.05,-0.95) -- (2.05,1) -- (1.95,1.1) -- (0,2.05);
			\draw (-1,-1) -- (-3,-1);
			\draw [fill=black] (-1.,-1.) circle (3pt);
			\draw[color=black] (-0.8,-0.8) node {$v_3$};
			\draw [fill=black] (1.,-1.) circle (3pt);
			\draw[color=black] (0.8,-0.8) node {$v_4$};
			\draw [fill=black] (-2.,1.) circle (3pt);
			\draw[color=black] (-1.7,0.9) node {$v_2$};
			\draw [fill=black] (2.,1.) circle (3pt);
			\draw[color=black] (1.7,0.925) node {$v_5$};
			\draw [fill=black] (0.,2.) circle (3pt);
			\draw[color=black] (0.025,1.75) node {$v_1$};
			\draw [fill=black] (-3.,-1.) circle (3pt);
			\draw[color=black] (-2.925,-0.7) node {$v^*$};
			\draw[color=black] (-0.9,1.35) node {$a$};
			\draw[color=black] (-1.3,0.1) node {$b$};
			\draw[color=black] (0,-0.8) node {$c$};
			\draw[color=black] (1.3,0.1) node {$d$};
			\draw[color=black] (0.9, 1.35) node {$e$};
			\draw[color=black] (-1.1,1.7) node {$e$};
			\draw[color=black] (-1.675,-0.075) node {$a$};
			\draw[color=black] (0,-1.25) node {$b$};
			\draw[color=black] (1.65,-0.075) node {$c$};
			\draw[color=black] (1.125, 1.75) node {$d$};
			\draw[color=black] (-2,-0.85) node {$a$};
		\end{tikzpicture}
		\begin{tikzpicture}
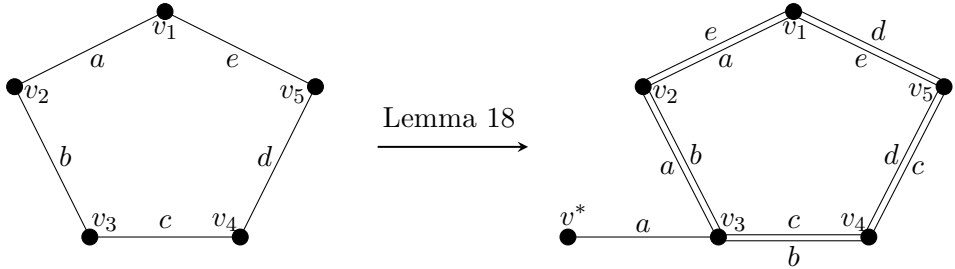
			\node at (0, -1.7) {\textbf{\hypertarget{figeight}{Figure 8}:} The ``growth'' of a rainbow $5$-cycle.};
			
		\end{tikzpicture}
	\end{center}
	
	\begin{proof}[Proof of \Cref{lem:5cyclegrowth}]
		Assume without loss that $\se_i \in \sD_i \in \cD$. Suppose $\se_i^+$ and $\se_i^-$ are the edges in $\sD_i$ satisfying $V(\se_i) \cap V(\se_i^+) = \{v_{i+1}\}$ and $V(\se_i) \cap V(\se_i^-) = \{v_i\}$, respectively. 
		
		Write $V \eqdef \{v_1, \dotsc, v_5\}$ for brevity. If there exists $v \in V(\se_i^{\bullet}) \setminus V$ for some $i \in [5]$ and $\bullet \in \{+, -\}$, say $i = 1$ and $\bullet = +$, then by choosing $v^* \eqdef v$ and $(j,k) = (0,1)$ the proof is done. 
		
		We thus assume that $V(\se_i^{\bullet}) \subseteq V$ for any $i \in [5]$ and $\bullet \in \{+, -\}$, and claim that this is impossible. To see this, we prove by contradiction. The following observation is quite useful: 
		
		\smallskip
		
		\hypertarget{fact}{\textbf{Fact.}} $V(\se_i^-) = \{v_{i-1}, v_i\}$ or $\{v_i, v_{i+3}\}$, and $V(\se_i^+) = \{v_{i+1}, v_{i+2}\}$ or $\{v_{i+1}, v_{i+3}\}$. 
		
		\emph{Proof of Fact}. Let $V(\se_i^-) \eqdef \{v_{i-1}, v'\}$. Then $v' \in \{v_{i-1}, v_{i+2}, v_{i+3}$\} is forced. However, $v' \neq v_{i+2}$, for otherwise $\se_i^-, \se_{i+2}, \se_{i+3}, \se_{i+4}$ form a rainbow $4$-cycle in $\cD$. The $V(\se_i^+)$ case is similar. $\blacksquare$
		
		\smallskip
		
		If $\se_i^+$ and $\se_{i+1}$ are coincident for all $i \in [5]$, then we cyclically shift the vertices via increasing $j$ by $1$ (note that the shift cannot happen indefinitely since the cycles $\sD_i$ are even). This does not change the situation, and so we may assume without loss that $V(\se_1^+) \neq \{v_2, v_3\}$. It follows from the \hyperlink{fact}{\underline{fact}} that $V(\se_1^+) = \{v_2, v_4\}$, which forces $V(\se^-_1) = \{v_1, v_5\}$, as shown in \hyperlink{fignineA}{Figure~9.A}. 
		
		We next look at $\se^{\pm}_2$. If $V(\se^-_2) = \{v_1,v_2\}$, then $\se_2^-, \se_1^+, \se_4, \se_5$ form a rainbow even cycle, a contradiction. So, $V(\se^-_2) = \{v_2,v_5\}$, and hence $V(\se^+_2) = \{v_3, v_4\}$ by the \hyperlink{fact}{\underline{fact}}, as shown in \hyperlink{fignineB}{Figure~9.B}. 
		
		We turn to $\se^{\pm}_3$ and $\se^{\pm}_5$ then. At this moment, we have a configuration that is symmetric in $(a, e)$ and $(b, c)$ (as seen in \hyperlink{fignineB}{Figure~9.B}). If $V(\se^-_3) = \{v_1, v_3\}$, then $\se_3^-, \se_2^+, \se_4, \se_5$ form a rainbow even cycle, a contradiction. So, the \hyperlink{fact}{\underline{fact}} implies $V(\se^-_3) = \{v_2,v_3\}$. By symmetry, $V(\se^+_5) = \{v_1, v_2\}$. We thus arrive at \hyperlink{fignineC}{Figure~9.C}. If $V(\se^+_3) = \{v_1, v_4\}$, then $\se_1, \se_2^-, \se_4, \se_3^+$ form a rainbow $4$-cycle, which is impossible. It then follows from the \hyperlink{fact}{\underline{fact}} and symmetry that $V(\se^+_3) = V(\se^-_5) = \{v_4, v_5\}$, as illustrated in \hyperlink{fignineD}{Figure~9.D}. 
		
		\begin{center}
			\begin{tikzpicture}[x=0.8cm,y=0.8cm]
				\clip(-2.5,-2.3) rectangle (2.5,2.5);
				\draw (0.,2.)-- (-2.,1.);
				\draw (-2.,1.)-- (-1.,-1.);
				\draw (-1.,-1.)-- (1.,-1.);
				\draw (1.,-1.)-- (2.,1.);
				\draw (2.,1.)-- (0.,2.);
				\draw (-2,1) -- (1,-1);
				\draw [fill=black] (-1.,-1.) circle (1.0pt);
				\draw[color=black] (-1,-1.2) node {$v_3$};
				\draw [fill=black] (1.,-1.) circle (1.0pt);
				\draw[color=black] (1,-1.2) node {$v_4$};
				\draw [fill=black] (-2.,1.) circle (1.0pt);
				\draw[color=black] (-2.25,1) node {$v_2$};
				\draw [fill=black] (2.,1.) circle (1.0pt);
				\draw[color=black] (2.25,1) node {$v_5$};
				\draw [fill=black] (0.,2.) circle (1.0pt);
				\draw[color=black] (0,2.2) node {$v_1$};
				\draw[color=black] (-1.05,1.65) node {$a$};
				\draw[color=black] (-1.75,-0.05) node {$\cancel{a}b$};
				\draw[color=black] (0,-1.15) node {$c$};
				\draw[color=black] (1.625,-0.025) node {$d$};
				\draw[color=black] (1.2, 1.65) node {$ae$};
				\draw[color=black] (-0.6,-0.125) node {$a$};
				\node at (0,-2) {\textbf{\hypertarget{fignineA}{Figure 9.A}}};
			\end{tikzpicture}
			\begin{tikzpicture}[x=0.8cm,y=0.8cm]
				\clip(-2.5,-2.3) rectangle (2.5,2.5);
				\draw (0.,2.)-- (-2.,1.);
				\draw (-2.,1.)-- (-1.,-1.);
				\draw (-1.,-1.)-- (1.,-1.);
				\draw (1.,-1.)-- (2.,1.);
				\draw (2.,1.)-- (0.,2.);
				\draw (-2,1) -- (1,-1);
				\draw (-2,1) -- (2,1);
				\draw [fill=black] (-1.,-1.) circle (1.0pt);
				\draw[color=black] (-1,-1.2) node {$v_3$};
				\draw [fill=black] (1.,-1.) circle (1.0pt);
				\draw[color=black] (1,-1.2) node {$v_4$};
				\draw [fill=black] (-2.,1.) circle (1.0pt);
				\draw[color=black] (-2.25,1) node {$v_2$};
				\draw [fill=black] (2.,1.) circle (1.0pt);
				\draw[color=black] (2.25,1) node {$v_5$};
				\draw [fill=black] (0.,2.) circle (1.0pt);
				\draw[color=black] (0,2.2) node {$v_1$};
				\draw[color=black] (-1.2,1.7) node {$a\cancel{b}$};
				\draw[color=black] (-1.75,-0.05) node {$\cancel{a}b$};
				\draw[color=black] (0,-1.2) node {$bc$};
				\draw[color=black] (1.625,-0.025) node {$d$};
				\draw[color=black] (1.2, 1.65) node {$ae$};
				\draw[color=black] (-0.6,-0.125) node {$a$};
				\draw[color=black] (0,1.2) node {$b$};
				\node at (0,-2) {\textbf{\hypertarget{fignineB}{Figure 9.B}}};
			\end{tikzpicture}
			\begin{tikzpicture}[x=0.8cm,y=0.8cm]
				\clip(-2.5,-2.3) rectangle (2.5,2.5);
				\draw (0.,2.)-- (-2.,1.);
				\draw (-2.,1.)-- (-1.,-1.);
				\draw (-1.,-1.)-- (1.,-1.);
				\draw (1.,-1.)-- (2.,1.);
				\draw (2.,1.)-- (0.,2.);
				\draw (-2,1) -- (1,-1);
				\draw (-2,1) -- (2,1);
				\draw[dash pattern=on 3pt off 3pt] (0,2) -- (-1,-1);
				\draw [fill=black] (-1.,-1.) circle (1.0pt);
				\draw[color=black] (-1,-1.2) node {$v_3$};
				\draw [fill=black] (1.,-1.) circle (1.0pt);
				\draw[color=black] (1,-1.2) node {$v_4$};
				\draw [fill=black] (-2.,1.) circle (1.0pt);
				\draw[color=black] (-2.25,1) node {$v_2$};
				\draw [fill=black] (2.,1.) circle (1.0pt);
				\draw[color=black] (2.25,1) node {$v_5$};
				\draw [fill=black] (0.,2.) circle (1.0pt);
				\draw[color=black] (0,2.2) node {$v_1$};
				\draw[color=black] (-1.2,1.65) node {$ae$};
				\draw[color=black] (-1.7,-0.025) node {$bc$};
				\draw[color=black] (0,-1.2) node {$bc$};
				\draw[color=black] (1.625,-0.025) node {$d$};
				\draw[color=black] (1.2, 1.65) node {$ae$};
				\draw[color=black] (-0.6,-0.125) node {$a$};
				\draw[color=black] (0,1.2) node {$b$};
				\draw[color=black] (-0.7,0.65) node {$\cancel{c}\cancel{e}$};
				\node at (0,-2) {\textbf{\hypertarget{fignineC}{Figure 9.C}}};
			\end{tikzpicture}
			\begin{tikzpicture}[x=0.8cm,y=0.8cm]
				\clip(-2.5,-2.3) rectangle (2.5,2.5);
				\draw (0.,2.)-- (-2.,1.);
				\draw (-2.,1.)-- (-1.,-1.);
				\draw (-1.,-1.)-- (1.,-1.);
				\draw (1.,-1.)-- (2.,1.);
				\draw (2.,1.)-- (0.,2.);
				\draw (-2,1) -- (1,-1);
				\draw (-2,1) -- (2,1);
				\draw[dash pattern=on 3pt off 3pt] (0,2) -- (1,-1);
				\draw[dash pattern=on 3pt off 3pt] (-1,-1) -- (2,1);
				\draw [fill=black] (-1.,-1.) circle (1.0pt);
				\draw[color=black] (-1,-1.2) node {$v_3$};
				\draw [fill=black] (1.,-1.) circle (1.0pt);
				\draw[color=black] (1,-1.2) node {$v_4$};
				\draw [fill=black] (-2.,1.) circle (1.0pt);
				\draw[color=black] (-2.25,1) node {$v_2$};
				\draw [fill=black] (2.,1.) circle (1.0pt);
				\draw[color=black] (2.25,1) node {$v_5$};
				\draw [fill=black] (0.,2.) circle (1.0pt);
				\draw[color=black] (0,2.2) node {$v_1$};
				\draw[color=black] (-1.2,1.65) node {$ae$};
				\draw[color=black] (-1.7,-0.025) node {$bc$};
				\draw[color=black] (0,-1.2) node {$bc$};
				\draw[color=black] (1.825,-0.025) node {$cde$};
				\draw[color=black] (1.2, 1.65) node {$ae$};
				\draw[color=black] (-0.6,-0.125) node {$a$};
				\draw[color=black] (0,1.2) node {$b$};
				\draw[color=black] (0.65,0.625) node {$\cancel{c}$};
				\draw[color=black] (0.45,-0.175) node {$\cancel{e}$};
				\node at (0,-2) {\textbf{\hypertarget{fignineD}{Figure 9.D}}};
			\end{tikzpicture}
		\end{center}
		
		Finally, we focus on $\se^{\pm}_4$. Indeed, we have $V(\se_4^-) = \{v_2, v_4\}$ or $\{v_3, v_4\}$ by the \hyperlink{fact}{\underline{fact}}. \hyperlink{fignineD}{Figure~9.D} shows that the former case generates a rainbow $4$-cycle on $\se_1, \se_4^-, \se_3^+, \se_5$ while the latter generates a rainbow $4$-cycle on $\se_2^-, \se_3^-, \se_4^-, \se_5^-$. We thus obtain the desired contradiction. 
	\end{proof}
	
	Assume $\fF_0$ and $\sf_0$ satisfy \Cref{lem:f'_5cycle}. Let $\sF_0$ be the forest part of $\fF_0$. That is, 
	\[
	\cP(\fF_0) = \{\sC_1, \dotsc, \sC_c, \sB_1, \dotsc, \sB_b, \sF_0\} \quad \text{with} \quad |\sF_0| = |\sF|. 
	\]
	From \Cref{lem:5cyclegrowth} we can find a subgraph of $\cD$ on six vertices $v_1, \dotsc, v_5$ and $v_*$. After some possible renaming of vertices, edges and colors, we assume this subgraph consists of the ingredients below:  
	\begin{itemize}
		\item $\widetilde{\sC} \eqdef \{\se_i = (v_iv_{i+1}, \alpha_i) : i \in [5]\}$ is the rainbow $5$-cycle in $\fF_0+\sf_0$ located by \Cref{lem:f'_5cycle}, and
		\item $\sp \eqdef (v^*v_1, \alpha_1)$ is a pendant edge of color $\alpha_1$ on vertices $v^*$ and $v_1$ located by \Cref{lem:5cyclegrowth}. 
	\end{itemize}
	
	We first claim that $\widetilde{\sC} - \sf_0 \subseteq \sF_0$. Since $\sf_0 \notin \fF_0$ and $\chi(\sf_0) \notin \chi(\fF_0)$, it follows from \Cref{lem:ectest} that $\widetilde{\sC}$ is edge-disjoint from $\sC_1, \dotsc, \sC_c$ and $\sB_1, \dotsc, \sB_b$. In particular, $\widetilde{\sC} - \sf_0 \subseteq \sF_0$. 
	
	We then claim that $\sp \notin \fF_0$. If $\sf_0 = \se_1$, then $\chi(\sp) = \chi(\sf_0) \notin \chi(\fF_0)$ follows from the choice of $\sf_0$ in \Cref{lem:f'_5cycle}, and so $\sp \notin \fF_0$. If $\sf_0 \in \{\se_2, \se_3, \se_4, \se_5\}$, then $\se_1 \in \sF_0$. This implies $\sp \notin \sF_0$ since $\sF_0$ is rainbow, and $\sp \notin \sC_i$, $\sp \notin \sB_j$ since $\sC_i, \sB_j$ are color-disjoint from $\sF_0$. We conclude that $\sp \notin \fF_0$. 
	
	Let $\widetilde{\sC} - \sf_0$ be a subgraph of $\sT_k \in \cP(\fF_0)$. Set $\fF_0' \eqdef \fF_0 + \sf_0 - \se_5$ and $\sF_0' \eqdef \sF_0 + \sf_0 - \se_5$. Note that $\sF_0'$ differs from $\sF_0$ only at $\sT_k'$, the rainbow tree from $\cP(\fF_0')$ containing $\widetilde{\sC}-\se_5$. Since $V(\sT_k') = V(\sT_k)$, 
	\[
	\cP(\fF_0') \eqdef \{\sC_1, \dotsc, \sC_c, \sB_1, \dotsc, \sB_b, \sF_0'\}
	\]
	shows that $\fF_0'$ is a Frankenstein subgraph of $\cD$ with $|\sF_0'| = |\sF_0| = |\sF|$. We remark that $\sp \notin \fF_0'$. 
	
	\begin{center}
		\begin{tikzpicture}[x=1.0cm,y=1.0cm]
			\clip(-6,-2.25) rectangle (4,5.5);
			\draw[color=blue] (0.,2.)-- (-2.,1.);
			\draw (-2.,1.)-- (-1.,-1.);
			\draw (-1.,-1.)-- (1.,-1.);
			\draw (1.,-1.)-- (2.,1.);
			\draw (2.,1.)-- (0.,2.);
			\draw[color=blue] (0.,4.)-- (0.,2.);
			\draw[color=blue] (-2,1) -- (-4,2) -- (-4,4) -- (-3,5) -- (-2, 5.25) -- (-1,5) -- (0,4);
			\draw[color=red] (-3,5) -- (0,2);
			\draw[color=red] (-4,4) arc[start angle=152,end angle=270,radius=3.4cm];
			\draw [fill=black] (-1.,-1.) circle (1.0pt);
			\draw[color=black] (-0.8,-0.8) node {$v_3$};
			\draw [fill=black] (1.,-1.) circle (1.0pt);
			\draw[color=black] (0.8,-0.8) node {$v_4$};
			\draw [fill=black] (-2.,1.) circle (1.0pt);
			\draw[color=black] (-1.25,0.9) node {$y_1 = v_2$};
			\draw [fill=black] (2.,1.) circle (1.0pt);
			\draw[color=black] (1.7,0.9) node {$v_5$};
			\draw [fill=black] (0.,2.) circle (1.0pt);
			\draw[color=black] (0.95,2.1) node {$y_{2\ell+2} = v_1$};
			\draw[color=blue] (-0.85,1.35) node {$\se_1$};
			\draw[color=blue] (-3.6,4.65) node {$\sq$};
			\draw[color=black] (0.95,1.35) node {$\se_5$};
			\draw[color=black] (-1.3,0.09) node {$\se_2$};
			\draw[color=black] (0,-0.8) node {$\se_3$};
			\draw[color=black] (1.3, 0.09) node {$\se_4$};
			\draw[color=black] (0,-1.25) node {$\widetilde{\sC}$};
			\draw [fill=black] (0.,4.) circle (1.0pt);
			\draw[color=black] (0.95,4.05) node {$y_{2\ell+1} = v^*$};
			\draw[color=black] (-1,5.2) node {$\cdots$};
			\draw[color=black] (-2,5.42) node {$\cdots$};
			\draw[color=black] (-4.2,2.15) node {$\vdots$};
			\draw[color=black] (-3,5.25) node {$y_{\mu}$};
			\draw[color=black] (-4.45,4) node {$y_{\mu-1}$};
			\draw[color=blue] (0.15,3) node {$\sp$};
			\draw [fill=black] (-4,2) circle (1.0pt);
			\draw [fill=black] (-4,4) circle (1.0pt);
			\draw [fill=black] (-3,5) circle (1.0pt);
			\draw [fill=black] (-1,5) circle (1.0pt);
			\draw [fill=black] (-2,5.25) circle (1.0pt);
			\draw[color=blue] (-1.975,4.925) node {$\sD_{\alpha_1}$};
			\draw[color=red] (-1.7,3.5) node {$\sP_1$};
			\draw[color=red] (-3.825,1) node {$\sP_2$};
			\draw[color=violet,dash pattern = on 3pt off 3pt] (-3.8,5.1) -- (0.8,0.5) -- (3,2.7);
			\draw[color=violet] (2.4,2.7) node {$\sS_1$};
			\draw[color=violet] (3,2.1) node {$\sS_2$};
			\node at (-1, -2) {\textbf{\hypertarget{figten}{Figure 10:}} An illustration of the proof of \Cref{thm:rainbow_ec}. };
		\end{tikzpicture}
	\end{center}
	
	Let $\sD_{\alpha_1}$ be the monochromatic even cycle from $\cD$ that contains $\se_1$ and $\sp$. Then 
	\[
	\sD_{\alpha_1} \eqdef \se_1 + (y_1y_2, \alpha_1) + \dotsb + (y_{2\ell+1}y_{2\ell+2}, \alpha_1) \qquad (y_1 \eqdef v_2, \, y_{2\ell+1} \eqdef v^*, \, y_{2\ell+2} \eqdef v_1, \, \ell \in \N_+). 
	\]
	By \Cref{lem:Fr_aux}, there are two connected components $\sS_1$ and $\sS_2$ of $\fF_0' - \se_1$ such that $v_1 \in V(\sS_1)$ and $v_2 \in V(\sS_2)$. Define~$\mu$ as the smallest index with $y_{\mu} \notin V(\sS_2)$. For similar reasons as ``$x_{\lambda-1} \in \sS_u$ and $x_{\lambda} \in \sS_w$'' in the proof of \Cref{lem:f'_5cycle}, we have that $y_{\mu-1} \in \sS_2$ and $y_{\mu} \in \sS_1$. Indeed, if $y_{\mu} \notin \sS_1$, then 
	\[
	\cP(\fF_0'+\se_5-\se_1+\sq) \eqdef \{\sC_1, \dotsc, \sC_c, \sB_1, \dotsc, \sB_b, \sF_0'+\se_5-\se_1+\sq\}, 
	\]
	where $\sq \eqdef (y_{\mu-1}y_{\mu}, \alpha_1)$, is another Frankenstein subgraph on $|\fF_*|+1$ edges, which contradicts \ref{max:edges}. 
	
	According to \Cref{prop:fpath}, we can find a rainbow path $\sP_1 \subseteq \sS_1$ with terminals $y_{\mu}$ and $v_1$. We can also find a rainbow path $\sP_2 \subseteq \sS_2$ whose terminals are $y_{\mu-1}$ and some $v_t \in \{v_2, v_3, v_4, v_5\}$. Assume further that $\sP_2$ is of minimum length. Note that $\sP_1 = \varnothing$ if $v_1 = y_{\mu}$, and $\sP_2 = \varnothing$ if $v_t = y_{\mu-1}$. 
	
	\hypertarget{claim}{\textbf{Claim.}} $\sP_1, \sP_2, \widetilde{\sC}$ are pairwise color-disjoint. 
	
	\emph{Proof of Claim.} Recall that $\widetilde{\sC}-\se_5 \subseteq \sT_k'$. If $\sP_1, \sP_2$ intersect some same part $\sG \in \cP(\fF_0')$, then it follows from the definitions of $\sS_1, \sS_2$ that $\sG = \sT_k'$. Since $\sT_k'$ is rainbow, and different parts in $\cP(\fF_0')$ are color-disjoint, we conclude that $\chi(\sP_1), \, \chi(\sP_2), \, \chi(\widetilde{\sC})$ are pairwise disjoint. $\blacksquare$
	
	\smallskip
	Decompose $\widetilde{\sC}$ into two rainbow paths $\widetilde{\sP}_1, \widetilde{\sP}_2$ with terminals $v_1$ and $v_t$ such that $\se_1 \in \widetilde{\sP}_1$. For instance, in \hyperlink{figten}{Figure~10} we have $v_t = v_3$, $\widetilde{\sP}_1 = \se_1+\se_2$ and $\widetilde{\sP}_2 = \se_3+\se_4+\se_5$. Set $\widetilde{\sP} \eqdef \sP_1 \cup \sP_2 \cup \{\sq\}$, and define $\widetilde{\sB} \eqdef \widetilde{\sP} \cup \widetilde{\sP}_1 \cup \widetilde{\sP}_2$. We are going to verify that $\widetilde{\sB}$ is a bad piece in three steps. 
	
	Firstly, we show that $\widetilde{\sB} = \widetilde{\sP} \cup \widetilde{\sP}_1 \cup \widetilde{\sP}_2$ is a theta graph with common terminals $v_1$ and $v_t$. Let $\widetilde{C}, P_1, P_2, q$ be the uncolored copies of $\widetilde{\sC}, \sP_1, \sP_2, \sq$, respectively. It suffices to show that $\widetilde{C}, P_1, P_2, \{q\}$ are pairwise disjoint. The definitions of $\sS_1, \sS_2$ indicate $q \notin P_1$, $q \notin P_2$ and $P_1 \cap P_2 = \varnothing$, $P_1 \cap \widetilde{C} = \varnothing$. The minimum-length assumption on $\sP_2$ implies $P_2 \cap \widetilde{C} = \varnothing$. To see that $q \notin \widetilde{C}$, we argue indirectly. If $q \in \widetilde{C}$, then $V(\sq) \subseteq V(\widetilde{\sC})$ and $y_{\mu} = v_1$. This implies $\sq=\sp$ hence $v^* \in \{v_1, \dotsc, v_5\}$, a contradiction. 
	
	Secondly, we prove that $\widetilde{\sP}$, $\widetilde{\sP}_1$, $\widetilde{\sP}_2$ are all rainbow, and that $\widetilde{\sB}$ is almost rainbow. Indeed, $\widetilde{\sP}_1, \widetilde{\sP}_2$ are rainbow because $\widetilde{\sC}$ is rainbow. Since $\chi(\sq) = \chi(\se_1) = \alpha_1$ and $\se_1 \in \widetilde{\sC}$, the \hyperlink{claim}{\underline{claim}} then implies that $\widetilde{\sP}$ is rainbow and $\widetilde{\sB}$ is almost rainbow. 
	
	Thirdly, we check that $|\widetilde{\sB}| \ge 7$. Since $\sq \in \widetilde{\sB}$, $\widetilde{\sC} \subseteq \widetilde{\sB}$ and $\sq \notin \widetilde{\sC}$, we obtain $|\widetilde{\sB}| \ge 6$. If $|\widetilde{\sB}| = 6$, then $V(\sq) \subseteq \{v_1, \dotsc, v_5\}$ and hence $y_{\mu-1} = v_t, \, y_{\mu} = v_1$, which contradicts $v^* \notin \{v_1, \dotsc, v_5\}$. So, $|\widetilde{\sB}| \ge 7$. 
	
	\smallskip
	If $|V(\widetilde{\sB}) \cap V(\sX)| \le 1$ for all $\sX \in \{\sC_1, \dotsc, \sC_c, \sB_1, \dotsc, \sB_b\}$, then the partition 
	\[
	\cP(\widetilde{\fF}) \eqdef \{\sC_1, \dotsc, \sC_c, \sB_1, \dotsc, \sB_b, \widetilde{\sB}\}
	\]
	exposes a \fw subgraph of $\cD$ with $c(\widetilde{\fF}) = c$ and $b(\widetilde{\fF}) > b$, which contradicts \ref{max:bpiece}. So, there exists $\sX_0 = \sC_i$ or $\sB_j$ such that $|V(\widetilde{\sB}) \cap V(\sX_0)| \ge 2$. 
	
	Consider the rainbow cycle $\widehat{\sC} \eqdef \widetilde{\sP} \cup \widetilde{\sP}_2$. We claim that $\chi(\sX_0) \cap \chi(\widehat{\sC} \setminus \sX_0) = \varnothing$. To see this, we begin with noticing that $\widehat{\sC}$ is a disjoint union $\sP_1 \cup \sP_2 \cup \widetilde{\sP}_2 \cup \{\sq\}$. The claim is then verified by: 
	\begin{itemize}
		\item $\chi(\sX_0) \cap \chi(\sP_1 \cup \sP_2 \setminus \sX_0) = \varnothing$ follows from $\sX_0 \in \cP(\fF_0')$ and $\sP_1 \cup \sP_2 \subseteq \fF_0'$; 
		\item $\chi(\sX_0) \cap \chi(\widetilde{\sP}_2 \setminus \sX_0) = \varnothing$ since $\sX_0 \in \cP(\fF_0), \, \chi(\sf_0) \notin \chi(\fF_0)$ and $\widetilde{\sP}_2 \subseteq \fF_0 + \sf_0$; 
		\item $\chi(\sq) \notin \chi(\sX_0)$ because $\chi(\sq) = \chi(\se_1) \in \chi(\sT_k')$ and $\sT_k', \sX_0$ are color-disjoint. 
	\end{itemize}
	It follows from $y_{\mu-1} \in \sS_2$ and $y_{\mu} \in \sS_1$ that $\sq \notin \fF_0'$, and hence $\sq \notin \sX_0$. Since $\sq \in \widetilde{P} \subseteq \widehat{\sC}$, we locate an edge $\sq \in \widehat{\sC} \setminus \sX_0$. \Cref{lem:ectest} then tells us that $|V(\widehat{\sC}) \cap V(\sX_0)| \le 1$, for otherwise a rainbow even cycle appears in $\cD$. Similarly, it follows from $\chi(\sX_0) \cap \chi(\widetilde{\sC} \setminus \sX_0) = \varnothing$, $\se_1 \in \widetilde{\sC} \setminus \sX_0$ and \Cref{lem:ectest} that $|V(\widetilde{\sC}) \cap V(\sX_0)| \le 1$. We thus obtain $|V(\widetilde{\sB}) \cap V(\sX_0)| = 2$ and $\widetilde{\sB} \cap \sX_0 = \varnothing$ by noticing $\widetilde{\sB} = \widehat{\sC} \cup \widetilde{\sC}$. 
	
	Suppose $V(\widetilde{\sB}) \cap V(\sX_0) \eqdef \{u, u_1\}$ with $u \in V(\widetilde{\sP}) \setminus V(\widetilde{\sP}_2)$ and $u_1 \in V(\widetilde{\sP}_1) \setminus V(\widetilde{\sP}_2)$. Denote by $\sP_{u, v_t}$ the subpath of $\widetilde{\sP}$ with terminals $u, v_t$, and by $\sP_{u_1, v_t}$ the subpath of $\widetilde{\sP}_1$ with terminals $u_1, v_t$. Write $\widehat{\sP} \eqdef \sP_{u, v_t} \cup \sP_{u_1, v_t}$. Then $\widehat{\sP}$ is a rainbow path because $\widehat{\sP} \subseteq \widetilde{\sB}$ and $\se_1 \notin \widehat{\sP}$. Since $\widetilde{\sB} \cap \sX_0 = \varnothing$, from \Cref{lem:ectest} we deduce that $\widehat{\sP} \cup \sX_0$ contains a rainbow even cycle, a contradiction. 
	
	\smallskip
	The proof of \Cref{thm:rainbow_ec} is complete. 
	
	\section{Concluding remarks}
	
	Write $\bracket{n} \eqdef \{3, 4, \dotsc, n\}$. For any positive integer $n$ and any $A \subseteq \bracket{n}$, let $f(n, A)$ be the minimum positive integer $N$ such that a rainbow $A$-cycle is guaranteed in every family of $N$ many $A$-cycles. It then follows from \Cref{thm:rainbow_c,thm:rainbow_oc,thm:rainbow_ec} that
	\[
	f(n, A) = \begin{cases}
		n \qquad &\text{when $A = \bracket{n}$}, \\
		2\lc\frac{n}{2}\rc-1 \qquad &\text{when $A = \bracket{n} \cap (2\Z+1)$}, \\
		\lf\frac{6(n-1)}{5}\rf+1 \qquad &\text{when $A = \bracket{n} \cap 2\Z$}. 
	\end{cases}
	\]
	We were unable to determine $f(n, A)$ when $A = \bracket{n} \cap (a\Z+b)$ in general. Another nice problem is to estimate $f(n, \{k\})$. It was proved independently by Gy\H{o}ri \cite{gyori} and Goorevitch, Holzman \cite{goorevitch_holzman} that $f(n, \{3\}) \approx \frac{n^2}{8}$. In particular, the value of $f(n, \{n\})$ concerning Hamiltonian cycles seems mysterious. 
	
	\section*{Acknowledgments}
	
	The first author is grateful to Boris Bukh, Ting-Wei Chao and Zilin Jiang for fruitful discussions. Part of this work was done after the graduation from Carnegie Mellon University in May 2023. As a postdoctoral researcher, the first author would like to thank
	\begin{itemize}
		\item Alfr\'{e}d R\'{e}nyi Institute of Mathematics (Budapest, Hungary) for hosting from Sep.~to Dec.~2023 (supported by ERC grant No.~882971, “GeoScape”, and the Erd\H{o}s Center), and
		\item Extremal Combinatorics and Probability Group (ECOPRO), Institute for Basic Science (IBS, Daejeon, South Korea) for hosting since Jan.~2024 (supported by IBS-R029-C4). 
	\end{itemize}
	
	The second author would like to thank Peking University for a pre-admission in his tenth grade, and to thank Beijing National Day School (high school) for allowing him to skip all regular classes in the academic year 2021--2022. These privileges resulted in plenty of free time to study all kinds of exciting new mathematics, especially to work on this problem on rainbow even cycles. 
	
	We thank two anonymous referees for their valuable feedback on earlier versions of this paper. 
	
	\bibliographystyle{plain}
	\bibliography{rainbow_ec}
	
\end{document}